\documentclass
[
    a4paper,
    DIV=11,
    abstracton,
]
{scrartcl}

\usepackage[utf8x]{inputenc}
\usepackage[T1]{fontenc}
\usepackage[english]{babel}
\usepackage[bf,normal]{caption}

\usepackage
{
    amsmath,
    amsfonts,
    amssymb,
    amsthm,
    mathtools,
    graphicx,
    xcolor,
    caption,
    subcaption,
    framed,
    authblk,
    enumitem,
    booktabs,
    kbordermatrix,
    todonotes
}

\usepackage[ruled,vlined,linesnumbered]{algorithm2e}

\usepackage[pdffitwindow=false,
            plainpages=false,
            pdfpagelabels=true,
            pdfpagemode=UseOutlines, 
            pdfpagelayout=OneColumn, 
            bookmarks=true,
            bookmarksopen=true,
            colorlinks=true,
            hyperfootnotes=false,
            linkcolor=blue,
            urlcolor=blue!30!black,
            citecolor=green!50!black]{hyperref}

\theoremstyle{plain} 
\newtheorem{thm}{Theorem}[section]

\theoremstyle{definition} 

\newtheorem{rem}[thm]{Remark}
\newtheorem{expl}[thm]{Example}
\newtheorem{textalgorithm}[thm]{Algorithm}

\newcommand{\R}{\mathbb{R}}										
\newcommand{\pd}[2]{\frac{\partial#1}{\partial#2}}              
\providecommand{\norm}[1]{\left\lVert #1 \right\rVert}          
\newcommand{\relmiddle}[1]{\mathrel{}\middle#1\mathrel{}}       

\newcommand\xqed[1]{\leavevmode\unskip\penalty9999 \hbox{}\nobreak\hfill \quad\hbox{#1}}
\newcommand{\exampleSymbol}{\xqed{$\triangle$}}

\newcommand{\ts}{\hspace*{0.1em}} 								

\allowdisplaybreaks

\author[1,2]{Jan-Hendrik Niemann}
\author[2,3]{Stefan Klus}
\author[1,2]{Christof Sch\"utte}
\affil[1]{Zuse Institute Berlin, Germany}
\affil[2]{Department of Mathematics and Computer Science, Freie Universit\"at Berlin, Germany}
\affil[3]{Department of Mathematics, University of Surrey, UK}
\title{Data-driven model reduction \\ of agent-based systems using \\ the Koopman generator}

\date{}

\begin{document}

\maketitle

\begin{abstract}
The dynamical behavior of social systems can be described by agent-based models. Although single agents follow easily explainable rules, complex time-evolving patterns emerge due to their interaction. The simulation and analysis of such agent-based models, however, is often prohibitively time-consuming if the number of agents is large. In this paper, we show how Koopman operator theory can be used to derive reduced models of agent-based systems using only simulation data. Our goal is to learn coarse-grained models and to represent the reduced dynamics by ordinary or stochastic differential equations. The new variables are, for instance, aggregated state variables of the agent-based model, modeling the collective behavior of larger groups or the entire population. Using benchmark problems with known coarse-grained models, we demonstrate that the obtained reduced systems are in good agreement with the analytical results, provided that the numbers of agents is sufficiently large.
\end{abstract}


\section{Introduction}

Systems of multiple agents that act and interact within a social network lead to complex dynamics and collective social phenomena. An agent can represent an individual person, a household, an organization, or any kind of discrete entity in an environment, which can be given, e.g., by geographical conditions, resources, infrastructure, but also rules or laws. Applications such as innovation spreading and infection kinetics (e.g., \cite{Kiesling2012, eubank+al2004}) range from data-based micro-simulations to abstract agent-based models (ABMs). A well-studied application concerns opinion dynamics and can be traced back to the \emph{voter model} introduced by Clifford and Sudbury \cite{Clifford1973}, developed in the 1970s. The name was coined by Holley and Liggett \cite{Holley1975} a few years later. In this model, an agent imitates the opinion of its neighbors. This means that whenever two agents with different opinions interact with each other, one of them copies the opinion of the other agent. There exist various modifications of the voter model, e.g., regarding the representation of the opinions, imitation, and interaction structure, see, for instance, \cite{jkedrzejewski2019statistical, redner2019, sirbu2017opinion, Herrerias-Azcue2019} for an overview.

Agent-based models provide an easily explainable and accessible framework for studying the dynamical behavior of interacting agents without requiring an extensive mathematical background. Models range from (highly detailed) microscopic stochastic descriptions following spatial movement and neighbor interactions \cite{conrad2018human} and individual-based stochastic descriptions in a network without movement \cite{Bolzern2019} to Markov chain approaches for collective population dynamics \cite{Banisch2012}. Most agent-based models have in common that they are hard to analyze due to their high-dimensionality. Additionally, simulations are often time-consuming so that a detailed analysis of such systems or parameter studies are typically infeasible. Especially for real-time decision and policy making this is clearly a disadvantage. One way to mitigate this is to compute surrogate models via machine learning approaches that can be used for calibration, sensitivity analysis, and parameter studies, see~\cite{Lamperti2018}. Another way is to represent the agents as a system of ordinary or stochastic (partial) differential equations (ODEs, mean-field ODEs, SDEs, or SPDEs), see, for instance, \cite{Kane2004, Kampen2007, grima2011accurate, Helfmann2021}. Assuming that the population of homogeneous agents that interact with each other (e.g., via a complete network) is sufficiently large, this system can be modeled as a Markov jump process (see also \cite{Banisch2012, Bolzern2019}), which in turn can be approximated using ordinary or stochastic differential equations~\cite{Herrerias-Azcue2019, Niemann2020}. This does not hold for all ABMs (consider, e.g., network-free or off-lattice models).

A drawback is that the aforementioned methods require knowledge about the process itself, which might not be available. Thus, there is a growing interest in learning the interaction laws of social dynamics in a data-driven fashion. One method is the so-called \emph{equation-free approach} pioneered by Kevrekidis et al.\ \cite{Kevrekidis2003, Kevrekidis2004}, which aims at circumventing the derivation of macroscopic, system-level equations when they are believed to exist but cannot be expressed in closed form. In~\cite{Zou2012}, the equation-free approach is used to obtain a reduced model of a spatio-temporally varying agent-based civil violence model. The obtained model is a stochastic differential equation that depends on two coarse-grained variables. The estimation of the drift and diffusion terms is accomplished by suitable short realizations of the agent-based simulation. Other applications of the equation-free approach are, e.g., bifurcation and stability analysis for ABMs or rare-event analysis \cite{Tsoumanis2010, Liu2015}. One key problem is the discovery of the right coarse-grained variables. If these are not known from physical insights or intuition, it is possible to use, e.g., a data-mining approach. In \cite{Liu2014}, the authors propose to use diffusion maps to learn the essential variables, resulting in an \emph{equation-free-variable-free} approach. In \cite{Lu2019}, a non-parametric approach for learning the interaction laws that is similar to parameter estimation problems for ordinary differential equations is proposed, assuming that the interaction depends only on pairwise distances between agents. Furthermore, it is shown that the learning rate is then independent of the dimension, making their approach suitable for large-scale systems. The data-driven approach described in \cite{WulkowKoltaiSchuette2020} utilizes memory terms to improve the accuracy of the coarse-grained model.

Our approach to learn coarse-grained systems for complex ABM dynamics relies on Koopman operator theory. The Koopman operator and its generator have been used for computing metastable and coherent sets, stability analysis, and control, but also for system identification, e.g., \cite{MauGon16, KKS16, Korda2018, Arbabi2018}. It was shown that by expressing the full-state observable in terms of the basis functions or eigenfunctions, it is possible to learn the governing equations of dynamical systems from data. While this has been mostly applied to ordinary differential equations \cite{Kaiser17, Kaiser18, MauGon16, MauGon17}, the approach can be naturally extended to stochastic differential equations, where the drift and diffusion terms are then estimated in a similar fashion \cite{Klus2020}. While Koopman operator-based methods have been successfully applied to molecular dynamics, fluid dynamics, engineering, and physics problems, the application of these methods to complex social systems such as ABMs, however, is still lacking, although notions like \emph{metastability} and \emph{coherence} exist in this context as well. The goal then is to study the coarse-grained behavior of complex ABMs based on data. If the model describes, for instance, the voting behavior of a large population, we are often not interested in each agent's decision but in the collective behavior of larger groups or the entire population. In~\cite{Fonoberova2018, Hogg2019}, the authors use Koopman mode analysis to investigate the dynamics of the spatial-temporal distribution of different agent types or to extract non-obvious information from the system's state indicating changes in the dynamics. Is was shown in \cite{Niemann2020} that the long-term characteristic behavior of ABMs can be determined by simulating (many) short trajectories of the corresponding SDE instead.

Our goal is to illustrate how coarse-grained models of complex ABM dynamics can be learned from data. The approach is based on~\cite{Klus2020}, with the difference that we here directly learn reduced models. Since we know the resulting limit processes in this case, which are given by a systems of ODEs or SDEs, we can compare the numerical results obtained for finitely many agents with the theoretical results. We demonstrate that under appropriate conditions the estimated models are in good agreement with known limit cases. The aim is to use the reduced models also for sensitivity analysis, parameter optimization, and control, by combining it with techniques proposed in \cite{Peitz2017, Korda2018, Peitz2020}. The main contributions of this work are:
\begin{itemize}
\item We show that the Koopman generator can be used to learn \emph{reduced} stochastic models from aggregated trajectory data that represents the collective behavior of larger groups or the entire population.
\item We demonstrate for a voter model defined on a complete network that the obtained reduced models are in good agreement with the SDE approximation for large population sizes and can not only be used for system identification but also for predictions of the temporal evolution. Furthermore, we show how the transition rate constants of the underlying Markov jump process corresponding to the ABM can be reconstructed.
\item We show that the proposed procedure also yields good reduced models that allow prediction in some other cases where the limit process is unknown or even far from a limit case. We demonstrate this for incomplete, clustered interaction networks (demonstrated again for the voter model) as well as models that do not have a network-based formulation (using a predator-prey model).
\end{itemize}
In general, this method requires a lot of data, which, however, is no problem in simulation studies where a surrogate model is required for the optimization or control of the full-complexity ABM.

The remainder of this paper is structured as follows: In Section~\ref{sec:Koopmanoperatortheory}, we introduce the stochastic Koopman operator, its generator, and \emph{generator extended dynamic mode decomposition} (gEDMD). We then briefly summarize the representation of ABMs as Markov jump processes and its SDE limit model for large population sizes in Section~\ref{sec:Modelingagentbasedsystems}. Furthermore, we introduce the voter model and the predator-prey model, which are used as guiding examples throughout the paper. In Section~\ref{sec:learningSDE}, we learn reduced models for complex ABM dynamics purely from aggregated data. We show in Section~\ref{sec:numericalresults} that, under certain conditions, the coarse-grained models agree with known limit cases. Furthermore, considering both ABMs with clustered interaction networks and ABMs without any underlying network structure, we demonstrate that the reduced models also allow prediction for other cases. Concluding remarks and future work will be discussed in Section~\ref{sec:conclusion}.

\section{Koopman operator theory}
\label{sec:Koopmanoperatortheory}

In this section, we will briefly introduce the stochastic Koopman operator, its generator, and generator EDMD, a variant of extended dynamic mode decomposition that can be used to approximate differential operators, see \cite{Klus2020} for details.

\subsection{The Koopman operator and its generator}

In what follows, let $\mathbb{X} \subset \R^d$ be the state space and $ f \in L^{\infty}(\mathbb{X}) $ a real-valued observable of the system, which can represent any kind of measurement. Furthermore, let $ \mathbb{E}[\,\cdot\,] $ denote the expected value. Given a stochastic differential equation of the form
\begin{equation} \label{eq:SDE}
	\mathrm{d}X_t = b(X_t) \ts \mathrm{d}t + \sigma(X_t) \ts \mathrm{d}W_t,
\end{equation}
where $ b \colon \R^d \to \R^d $ is the drift term, $ \sigma \colon \R^d \to \R^{d \times s} $ the diffusion term, and $ W_t $ an $ s $-dimensional Wiener process, the stochastic Koopman operator is defined by
\begin{equation*}
	(\mathcal{K}^t f)(x) = \mathbb{E}[f(\Phi^t(x))].
\end{equation*}
Here, $ \Phi^t $ is the flow map associated with \eqref{eq:SDE}. It can be shown that the infinitesimal generator of the stochastic Koopman operator is
\begin{equation*}
	\mathcal{L} f
	= \sum_{i=1}^d b_i \ts \pd{f}{x_i} + \frac{1}{2} \sum_{i=1}^d \sum_{j=1}^d a_{ij} \ts \pd{^2 f}{x_i \ts \partial x_j},
\end{equation*}
where $a = \sigma \ts \sigma^\top$. The adjoint operator is given by
\begin{equation*}
	\mathcal{L}^* f = -\sum_{i=1}^d \pd{(b_i \ts f)}{x_i}  + \frac{1}{2} \sum_{i=1}^d \sum_{j=1}^d \pd{^2 (a_{ij} \ts f)}{x_i \ts \partial x_j}.
\end{equation*}
The function $ u(t, x) = \mathcal{K}^t f(x) $ solves the \emph{Kolmogorov backward equation} given by the second-order partial differential equation $ \pd{u}{t} = \mathcal{L} u $, see \cite{Met07}. Moreover, $ \pd{u}{t} = \mathcal{L}^* u $ is called \emph{Fokker--Planck equation}~\cite{LaMa94}. For deterministic dynamical systems, $ \sigma \equiv 0 $ and consequently also $ a \equiv 0 $ so that we obtain a first-order partial differential equation, namely the Liouville equation.

\subsection{Infinitesimal generator EDMD}
\label{sec:gEDMD}

While the classical \emph{extended dynamic mode decomposition} (EDMD) approximates the Koopman operator or the Perron--Frobenius operator~\cite{WKR15, KKS16}, we now seek to approximate their generators from data. We thus introduce \emph{generator EDMD} or, in short, gEDMD, which was proposed in~\cite{Klus2020}. Assume that we have $m$ measurements of the system's state $ \{\ts x_l \ts\}_{l=1}^m $, its drift $ \{\ts b(x_l) \ts\}_{l=1}^m $, and diffusion $ \{\ts \sigma(x_l) \ts\}_{l=1}^m $. We will discuss in Section~\ref{sec:learningSDE} how to obtain these pointwise estimates. Then, choosing a set of basis functions $ \{\ts \psi_i \ts\}_{i=1}^n $, which is sometimes also called \emph{dictionary}, and writing it in vector form as $ \psi(x) = [\psi_1(x), \dots, \psi_n(x)]^\top $, we define
\begin{equation*}
\mathrm{d}\psi_k(x) = (\mathcal{L} \psi_k)(x) = \sum_{i=1}^d b_i(x) \ts \pd{\psi_k}{x_i}(x) + \frac{1}{2} \sum_{i=1}^d \sum_{j=1}^d a_{ij}(x) \ts \pd{^2 \psi_k}{x_i \ts \partial x_j}(x).
\end{equation*}
For all measurements and basis functions, we can now assemble the matrices
\begin{equation}\label{eq:gEDMDmatrices}
	\Psi_X =
	\begin{bmatrix}
	\psi_1(x_1) & \dots  & \psi_1(x_m) \\
	\vdots      & \ddots & \vdots      \\
	\psi_n(x_1) & \dots  & \psi_n(x_m)
	\end{bmatrix}
	\quad \text{and} \quad
	\mathrm{d}\Psi_X =
	\begin{bmatrix}
	\mathrm{d}\psi_1(x_1) & \dots  & \mathrm{d}\psi_1(x_m) \\
	\vdots                & \ddots & \vdots                \\
	\mathrm{d}\psi_k(x_1) & \dots  & \mathrm{d}\psi_k(x_m)
	\end{bmatrix},
\end{equation}
where $ \Psi_X, \mathrm{d}\Psi_X \in \R^{n \times m} $. Assuming there exists a matrix $M$ such that $ \mathrm{d}\Psi_X = M \Psi_X $, we solve the problem in the least-square sense by minimizing $\norm{\smash{\mathrm{d}\Psi_X - M \Psi_X}}_F $ since in general this problem cannot be solved exactly. Here, $\norm{\cdot}_F$ denotes the Frobenius norm. The least-squares solution is given by
\begin{equation*}
M = \mathrm{d}\Psi_X \Psi_X^+ = \big(\mathrm{d}\Psi_X \Psi_X^\top\big) \big(\Psi_X \Psi_X^\top\big)^+,
\end{equation*}
where $A^+$ denotes the Moore--Penrose pseudoinverse of a matrix $A$. The matrix $ L = M^\top $ is an empirical estimate of the matrix representation of the infinitesimal generator $\mathcal{L}$ as shown in \cite{Klus2020}. In the infinite data limit, gEDMD converges to a Galerkin approximation of the generator, i.e., a projection onto the space spanned by the basis functions.

\subsection{System identification}
\label{sec:systemidentification}

Let $\mathbb{X}$ be bounded so that the full-state observable $g(x) = x$ is (component-wise) contained in $L^{\infty}(\mathbb{X})$.
With the aid of the full-state observable, it is possible to reconstruct the governing equations of the underlying dynamical system. We assume that the function $ g(x) = x $ can be represented by the basis functions $ \psi $. The easiest way to accomplish this is to add the observables $ \{ \ts x_i \ts \}_{i=1}^d $ to the dictionary. Let $ B \in \mathbb{R}^{n \times d} $ be the matrix such that $ g(x) = B^\top \ts \psi(x) $. The system can directly be represented in terms of the basis functions,
\begin{equation*}
	(\mathcal{L} g)(x) = b(x) \approx (L B)^\top \ts \psi(x),
\end{equation*}	
which, for a deterministic dynamical system, is equivalent to SINDy~\cite{Brunton2015}. For non-deterministic systems and for $\psi_k(x) = x_ix_j$, note that the diffusion term can be identified by
\begin{equation}\label{eq:diff_identification}
	a_{ij}(x) \approx (\mathcal{L}\psi_k)(x) - b_i(x) x_j - b_j(x) x_i,
\end{equation}
provided that $b_i$ and $b_j$ as well as $b_i(x)x_j$ and $b_j(x)x_i$ are contained in the space spanned by the basis functions. If the drift term $\sigma$ itself is needed, we can obtain it using a Cholesky decomposition of $a$, see \cite{Klus2020}.

\section{Modeling agent-based systems}
\label{sec:Modelingagentbasedsystems}

We consider agent-based systems of $N$ interacting agents. For each system, there is a set $\{S_1, \dots, S_d\}$ of types available to the agents, a set $\{R_1, \dots, R_K\}$ of transition rules that define possible changes between the types $S_i$, and a set of propensity functions specifying the rates of random occurrences of the transitions. The ABM state space is given by $\{1,\dots,d\}^N$ and grows like $d^N$, which is problematic for large $N$. For this reason, we describe the ABM via the population state, i.e., we count the number of agents of each type. The population state space grows like $N^d$ in the worst case. If we assume random interactions between all agents (e.g., via a complete network) and indistinguishable agents, then the population state space description is exact. In all other cases it involves an approximation error due to aggregation of the ABM state space. 

We will consider two different agent-based models and modeling approaches. The first one is a continuous-time voter model without spatial resolution where the agents are nodes in an interaction network and each of them can switch between $d$ different types according to some given transition rules. This model is similar to the discrete-time model in \cite{Banisch2012}. The second ABM is a spatial (i.e., there is no underlying network) predator-prey model formulated in discrete-time. Unlike in the first model, the agents are not changing their types (in this context called \emph{breed}). Instead, transitions in the population state are caused by reproduction and death of predators and prey. The population size is thus not constant.

We will now describe the representation of agent-based systems (using the population state) as a Markov jump processes and their approximation by SDEs for large population sizes. For further details, we refer the reader to \cite{Niemann2020}.

\subsection{Agent-based models as Markov jump processes}
\label{sec:ABMasMJP}

At any time $t$, the population state $x \in \mathbb{X}$ of the ABM is fully described by the vector
\begin{equation*}
	x = [x_1,\dots, x_d]^\top \in \mathbb{N}_0^d,
\end{equation*}
where $x_i$ is the number of agents of type $S_i$. For the sake of simplicity, we assume in this subsection random interactions between all agents so that transitions between agent types imply transitions between population states. We use a formalism that is most commonly used in the chemical context, where each transition rule is represented by an equation of the form
\begin{equation*}
	R_k: ~~ a_{1k} \ts S_1 + \dotsc + a_{dk} \ts S_d ~~\mapsto~~  b_{1k} \ts S_1 + \dotsc + b_{dk} \ts S_d.
\end{equation*}
It induces an instantaneous change in the system's state of the form $x \mapsto x + \nu_k$, where $\nu_k = (\nu_{ik})_{i = 1,\dots,d }$, defined by $\nu_{ik} \coloneqq b_{ik} - a_{ik}$, describes the net change in the number of agents of each type $S_i$ due to transitions $R_k$. Transition $R_k$ occurs in an infinitesimal time step $\mathrm{d}t$ with probability $\alpha_{k}(x)\ts\mathrm{d}t$, where $\alpha_{k} \colon \mathbb{X} \to [0, \infty)$ denotes the \emph{propensity function} associated with transition $R_k$. We assume that the propensity $\alpha_k$ is proportional to the number of combinations of agents in $x$, and, moreover, that it scales with the total population size $N$, i.e., 
\begin{equation*}
	\alpha_k(x) = 
	\begin{cases} \displaystyle
		\gamma_k N\prod_{i=1}^d  \frac{1}{N^{a_{ik}}} \binom{x_i}{a_{ik}} , & \text{if } x_i \geq a_{ik} \text{ for all } i=1,\dots,d, \\[3ex]
		0, & \text{otherwise}.
	\end{cases}
\end{equation*}
Here, $\gamma_k > 0$ denotes the rate constant for the $k$th transition $R_k$.

The evolution of the population state can be described by a continuous-time stochastic process $\{X_t\}_{t \ge 0}$ with
\begin{equation*}
	X_t = (x_i(t))_{i=1,\dots,d} \in \mathbb{X},
\end{equation*}
where $x_i(t)$ denotes the number of agents of type $S_i$ at time $t$. It is a Markov jump process, i.e., it is piece-wise constant with jumps of the form $X_t \mapsto X_t + \nu_{k}$.

Let $P(x,t) \coloneqq \mathbb{P}[X_t=x\mid X_0=x_0]$ denote the probability of finding the process in state $x$ at time $t$ given some initial state $x_0$. The temporal evolution of $\{X_t\}_{t \ge 0}$ can then be described by the \emph{Kolmogorov forward equation} given by
\begin{equation} \label{eq:CME}
	\frac{\mathrm{d}P(x,t)}{\mathrm{d} t}  =  \sum_{k=1}^K \big[\alpha_{k}(x-\nu_k)P(x-\nu_k,t)-\alpha_{k}(x)P(x,t)\big].
\end{equation}
By setting $\alpha_{k}(x) \coloneqq 0$ and $P(x,t) \coloneqq 0$ for $x\notin \mathbb{N}_0^d$, we exclude terms in the right-hand side of \eqref{eq:CME} where the argument $x - \nu_k$ contains negative entries. Since in general the Kolmogorov forward equation of the ABM process cannot be solved analytically, the distribution of the process can be estimated by Monte Carlo simulations, which can be generated using Gillespie's stochastic simulation algorithm \cite{gillespie1976general}.

Assuming convergence of the propensity functions for $N\to\infty$, it is well-known that the rescaled jump process $X_tN^{-1}$ converges to the frequency process $C(t)$, $t\ge 0$, given by the SDE
\begin{equation}\label{eq:SDE_limit}
	\mathrm{d}C (t) = \sum_{k=1}^K \widetilde{\alpha}_k(C(t)) \ts \nu_{k} \ts \mathrm{d}t + \sum_{k=1}^K \frac{1}{\sqrt{N}} \sqrt{\widetilde{\alpha}_k(C(t))} \ts \mathrm{d}W_k(t) \ts \nu_{k},
\end{equation}
with initial state $C(0) = \lim_{N \to \infty} X_0 N^{-1}$, independent Wiener processes $W_k(t)$, $k=1,\dots,K$, and rescaled propensities, i.e., $\widetilde{\alpha}_k(c) = N^{-1} \alpha_{k}(cN)$ \cite{Kurtz1976}. The SDE limit model \eqref{eq:SDE_limit} is also known as the \emph{chemical Langevin equation} in the context of chemical reaction kinetics \cite{gillespie2000chemical}. Written as an SDE of the form \eqref{eq:SDE}, the drift and diffusion terms $b(c)$ and $\sigma(c)$ are given by
\begin{align}
	b(c) & = \sum_{k=1}^K \widetilde{\alpha}_k(c) \ts \nu_{k}, \label{eq:SDEdrift} \\
	\sigma(c) & = \frac{1}{\sqrt{N}}
	\begin{bmatrix}
		\sqrt{\widetilde{\alpha}_1(c)} \ts \nu_{1} & \dots & \sqrt{\widetilde{\alpha}_K(c)} \ts \nu_{K}
	\end{bmatrix}. \label{eq:SDEdiffusion}
\end{align}

\subsection{Extended voter model}
\label{exp:guiding}

\begin{figure}[!t]
	\centering
	\begin{subfigure}[t]{0.45\textwidth}
		\centering
		\caption{}
		\includegraphics[scale=0.9]{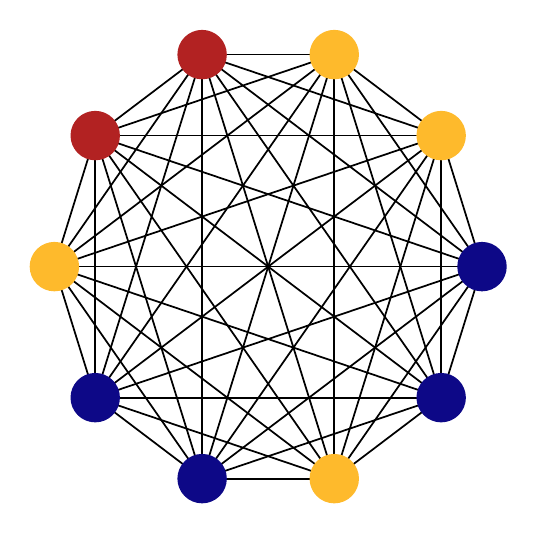}
		\vspace{0.75cm}
	\end{subfigure}
	\hfill
	\begin{subfigure}[t]{0.45\textwidth}
		\centering
		\caption{}
		\includegraphics[scale=0.9]{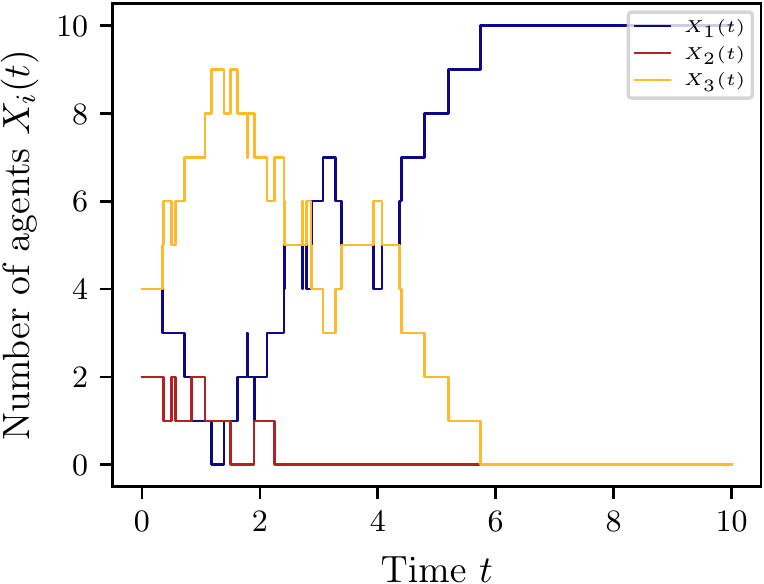}
	\end{subfigure}
	\caption{(a) Network of $N=10$ agents having three different types (blue, red yellow) and (b)~a possible trajectory of the jump process for the rate constants $\gamma_{12} = \gamma_{23} = \gamma_{31} = 2$, $\gamma_{32} = \gamma_{21} = \gamma_{13} = 1$ and $\gamma_{ij}' = 0.01$ for $i,j = 1,\dots,3$.}
	\label{fig:guiding}
\end{figure}

Throughout the paper, we will use the \emph{extended voter model} (EVM) with $N$ agents, $d$ types, and two sorts of transitions as one of two guiding examples. This model is well-known, e.g., as the noisy multi-state voter model, for describing foraging ant colonies, or chemical systems, see \cite{Herrerias-Azcue2019, Biancalani2014, Ohkubo2008}. The agents are the nodes in an interaction network. Given two agents with types $S_i \neq S_j$, \emph{imitation} or \emph{adaption} is a second-order transition of the form $R_{ij} \colon S_i + S_j \mapsto 2 \ts S_j$, whereas \emph{exploration} or \emph{mutation} is a first-order transition of the form $R_{ij}' \colon S_i \mapsto S_j$. Imitation happens whenever one agents of type $S_i$ adopts the type of another agent with different type $S_j$. It can be interpreted as adopting an opinion or technology, or also as being infected. Exploration corresponds to an independent change of the agent's type. Given a complete network, the propensity functions for imitative and exploratory transitions $R_{ij}$ and $R_{ij}'$ are given by
\begin{equation*}
	\alpha_{ij}(x) = \frac{\gamma_{ij}}{N} x_ix_j \quad
	\text{and} \quad \alpha_{ij}' = \gamma'_{ij} x_i,
\end{equation*}
where $\gamma_{ij}, \gamma'_{ij} > 0$ denote the rate constants for the transitions. Figure~\ref{fig:guiding}~(a) shows a graph with $N=10$ nodes representing the interaction network. Here, the agents can have three different types (represented by blue, red, and yellow vertices). Figure~\ref{fig:guiding}~(b) shows a trajectory of the Markov jump process.

\subsection{Predator-prey model}
\label{sec:PPM}

The second agent-based model we consider as a guiding example in this work is a predator-prey model (PPM), where the agents move freely in a given domain. We formulate the PPM in discrete time using intuitive text-based transition rules to emphasize its ABM character.

\begin{figure}[!t]
	\centering
	\includegraphics[scale=1]{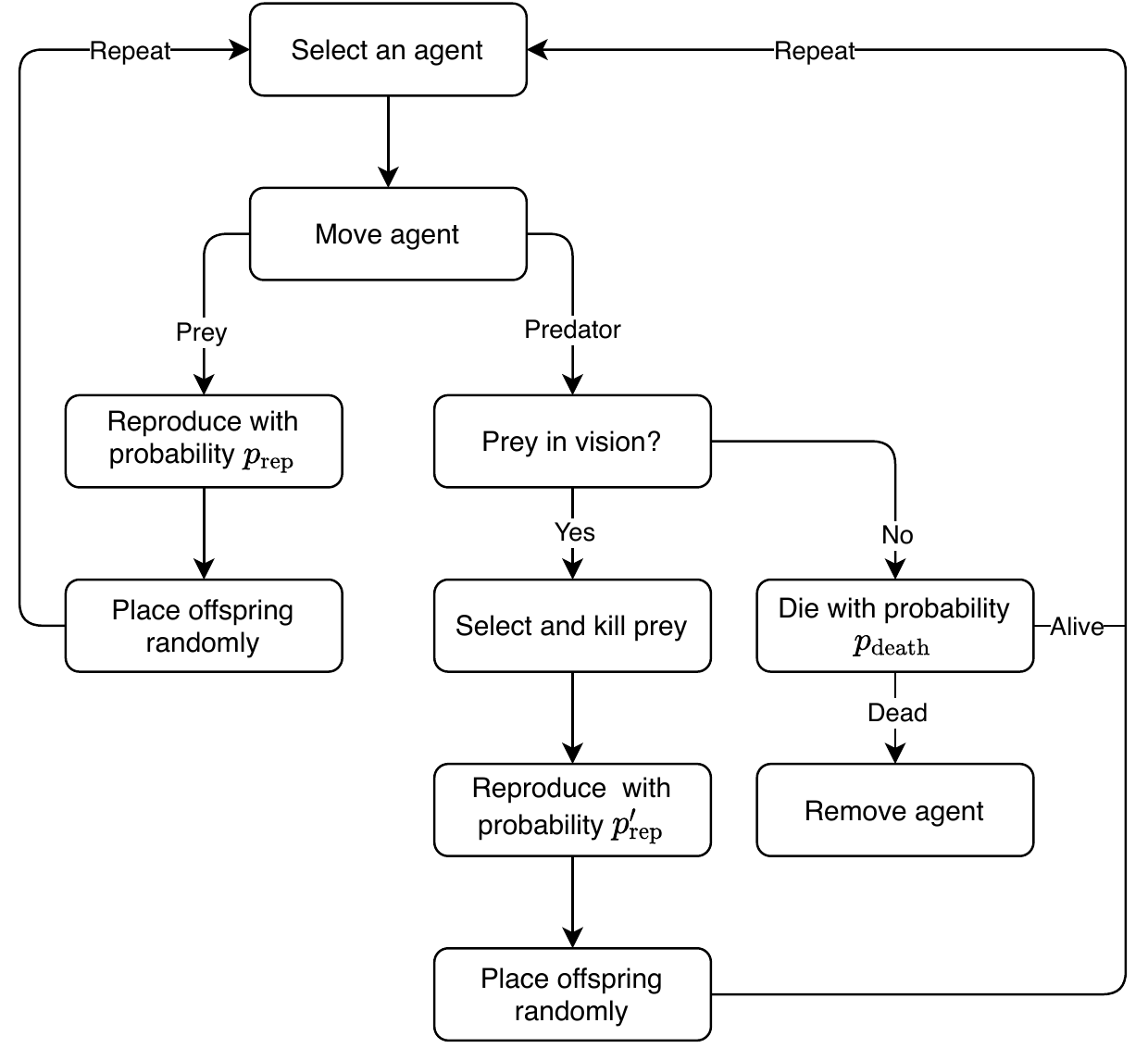}
	\caption{Flow chart of the predator-prey model.}
	\label{fig:PPMflowchart}
\end{figure}

\begin{figure}[!t]
	\centering
	\begin{subfigure}[t]{0.45\textwidth}
		\centering
		\caption{}
		\includegraphics[scale=0.9]{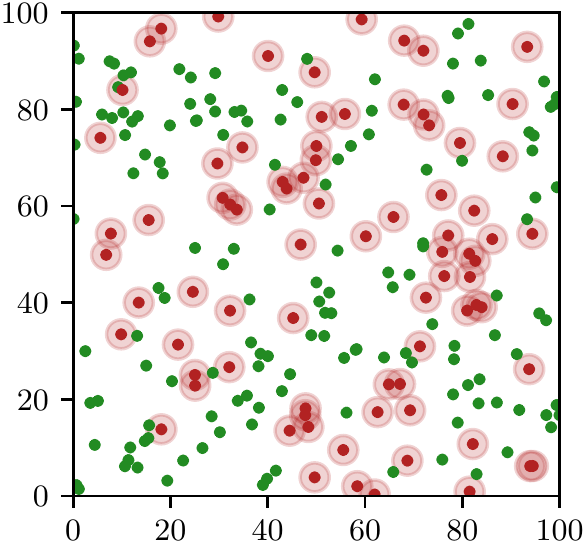}
		\vspace{0.5cm}
	\end{subfigure}
	\hfill
	\begin{subfigure}[t]{0.45\textwidth}
		\centering
		\caption{}
		\includegraphics[scale=0.9]{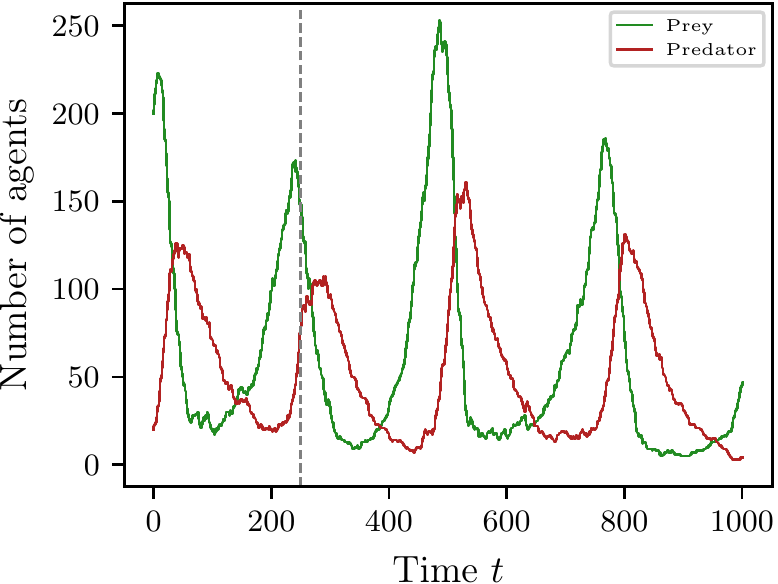}
	\end{subfigure}
	\caption{(a) Snapshot of the state of the predator-prey ABM at time $t=250$. Red and green dots represent predators and prey, respectively. The radius of vision is indicated by the light-red area around the predators. (b)~Simulation of the predator-prey model for the parameters given in Table \ref{tab:PPM_parameters} on page~\pageref{tab:PPM_parameters}. The vertical gray dashed line indicates the time where the snapshot in (a) is taken.}
	\label{fig:PPM}
\end{figure}

Given a continuous and periodic space, all agents are constantly performing Gaussian random walks with normally distributed step size. This means, given its current position in space $x_i(k_0)$, after $k$ time steps the agent is located at position $x_i(k_0+k) = x_i(k_0) + \sum_{i=1}^{k} \xi_i$ for $\xi_i \sim \sqrt{h}\ts\ts\mathcal{N}(0, 1)$. There are two breeds of agents: predator agents and prey agents. We will denote them as \emph{predators} and \emph{prey}, respectively. At each time step, all agents carry out the following steps corresponding to their breed:
\begin{itemize}
	\item A prey moves and reproduces with probability $p_\text{rep}$. The offspring is placed randomly in the space. 
	\item A predator moves and looks for prey within a radius of vision $v$. If there is prey within the radius of vision, the predator chooses its victim randomly and kills it. The predator can only reproduce with probability $p_\text{rep}'$ if it killed a prey before. The offspring is placed randomly in the space. If there is no prey in the radius of vision, the predator dies with probability $p_\text{death}$.
\end{itemize}
A flow chart describing the PPM in more detail can be found in Figure~\ref{fig:PPMflowchart}. In the absence of predators, the prey has an unlimited growth, which can be interpreted as independence of resources. There is no competition between the prey. The growth is only kept in check by the existence of predators. The population size is clearly not constant here. Figure~\ref{fig:PPM}~(a) shows a snapshot of the PPM for a realization using the parameters summarized in Table \ref{tab:PPM_parameters}. Green and red dots represent prey and predators, respectively. The search radius for prey is indicated by the light-red area around the red dots. The aggregate state is given by the number of prey and predators, respectively.

\begin{rem}\label{rem:PPMasMJP}
	Due to the spatial component of the PPM, it cannot be formulated directly using the formalism summarized in Section~\ref{sec:ABMasMJP}. Assuming a \emph{well-mixed} system and denoting prey by $S_1$ and predators by $S_2$, the rules given above translate to
    \begin{equation*}
		\begin{array}{rcll}
			S_1 & \mapsto & 2 \ts S_1, & \quad\text{(reproduction of prey)} \\
			S_1 + S_2 & \mapsto & 2 \ts S_2, & \quad\text{(reproduction of predators)} \\
			S_2 & \mapsto & \emptyset, & \quad\text{(death of predators)}
		\end{array}
    \end{equation*}
    for some rate constants $\gamma_i > 0$, $i = 1,\dots, 3$. Then the aggregate state of the PPM resembles the stochastic Lotka--Volterra predator-prey differential equations.
\end{rem}

In the next section, we will show how we can obtain reduced models of agent-based models using simulation data only.

\section{Learning coarse-grained models from data}
\label{sec:learningSDE}

We will now illustrate how to learn reduced models for large agent-based dynamics from aggregated trajectory data using the Koopman generator. The approach is based on \cite{Klus2020}. First, we estimate drift and diffusion \emph{pointwise}, cf.~Section~\ref{sec:gEDMD}. Subsequently, we apply gEDMD to the estimates to obtain a \emph{global} description of the drift and diffusion terms, cf.~Section~\ref{sec:systemidentification}. For the EVM, we will show in Section~\ref{sec:complete_networks} that the identified SDE coincides with the SDE limit model \eqref{eq:SDE_limit}, provided that the number of agents is sufficiently large. We will now go through the main steps that are necessary to learn the Koopman generator from data generated by an ABM.

\paragraph{Measurements.} Assume that we have access to $m$ measurements of an aggregate state variable of a given ABM. This aggregate state can represent the number of agents sharing, e.g., the same type $S_i$ or belonging to some group. These $m$ measurements will be the starting point. Let us denote them by $\{x_l\}_{l=1}^m$. If possible, we choose the measurements $x_l$ such that they are uniformly distributed in the aggregate state space $\mathbb{X}$ to ensure a good coverage of the whole (aggregate) space. One way to achieve this is by constructing an appropriate map from the macroscopic (aggregate) state to the microscopic ABM state. By \emph{appropriate} we mean that the mapped macroscopic state and a naturally developed ABM state with same aggregate variables agree in probability. Practically, this means that if, e.g., the agents follow a certain spatial distribution, this needs to be taken into account when constructing the map. Another, rather straightforward, possibility is to gather the measurements ``on the fly'', i.e., by using the states belonging to trajectories obtained from the simulation of the ABM.

\paragraph{Pointwise estimates.} Since the drift and diffusion terms $b$ and $\sigma$ are in general unknown, we estimate them pointwise via finite difference approximations for each measurement $ \{\ts x_l \ts\}_{l=1}^m $ using the Kramers--Moyal formulae
\begin{subequations}\label{eq:kramersmoyal}
	\begin{align}
		b(x) &\coloneqq \lim_{\tau \to 0} \mathbb{E}\left[ \frac{1}{\tau} (X_\tau - x) \relmiddle| X_0 = x \right], \\
		a(x) &\coloneqq \lim_{\tau \to 0} \mathbb{E}\left[ \frac{1}{\tau} (X_\tau - x)(X_\tau - x)^\top \relmiddle| X_0 = x \right].
	\end{align}
\end{subequations}
The formulae can be deduced from the Kramers--Moyal expansion, see, e.g.,~\cite{Risken1996}. These expressions can be evaluated by Monte Carlo methods via multiple short trajectories at each data point $ \{\ts x_l \ts\}_{l=1}^m $. The simulation of multiple short realizations of the original ABM is comparable to the equation-free approach and common practice in the context of transfer operator approximations. These pointwise estimates of the drift and diffusion for each training data point form the first stage to obtain a global description of them via gEDMD.

\paragraph{Conservation laws.} If the aggregate state is subject to a conservation law, e.g., if the number of agents is constant for all time $t\ge 0$, we have only $d-1$ degrees of freedom and the aggregated trajectory data belongs to a $d-1$ dimensional system, i.e., the number of agents $x_j(t)$ can be expressed by
\begin{equation*}
	x_j(t) = N - \sum_{i \neq j} x_i(t).
\end{equation*}
We thus reduce each measurement by keeping, without loss of generality, only the first $d-1$ entries. This eliminates redundant representations of the system. Additionally, we can scale the measurements by the number of agents, $N$, to obtain a frequency representation $c_i(t) = \frac{x_i(t)}{N}$. 

\paragraph{Basis functions.}Next, we need to choose a set of basis functions $\{ \psi_i\}_{i=1}^n$. This is a non-trivial step since in general it is not clear how the drift term $b$ and diffusion term $\sigma$ of the SDE \eqref{eq:SDE} look like. If we assume that the SDE approximation of the ABM adheres to the model structure introduced in Section~\ref{sec:Modelingagentbasedsystems} and comprises at most $p$th order transitions, we can show that monomials of degree up to $p+1$ are sufficient to correctly identify the model of the form~\eqref{eq:SDE_limit}. The highest order transition coincides with the maximum degree of all propensity functions. First, to identify the drift term \eqref{eq:SDEdrift}, we conclude from the propensity functions that the set of basis functions needs to contain at least monomials up to degree $p$. Second, as gEDMD identifies $a = \sigma \ts \sigma^\top$ and not the diffusion term \eqref{eq:SDEdiffusion} itself, we obtain for $c=x/N$
\begin{align*}
	a(c) \coloneqq \sigma(c) \ts \sigma(c)^\top = \sum_{k=1}^K \frac{1}{N} \ts \widetilde{\alpha}_k(c) \ts \nu_{k}\nu_{k}^\top,
\end{align*}
which shows that monomials are sufficient for the identification of the diffusion term as well. Finally, to identify the diffusion term via \eqref{eq:diff_identification}, we argue that also monomials of degree $p+1$ are needed.

\paragraph{Identification.}We are now able to assemble the matrices $\Psi_X$ and $\mathrm{d}\Psi_X$ in \eqref{eq:gEDMDmatrices} and solve the minimization problem $\norm{\smash{\mathrm{d}\Psi_X - M \Psi_X}}_F $ to obtain an approximation $L = M^\top$ of the infinitesimal generator $\mathcal{L}$ associated with the ABM. For a suitable projection matrix $B$, we identify the drift and diffusion terms. These are now global descriptions (i.e., functions depending on $x$) forming the second stage, cf. Section~\ref{sec:systemidentification}. The overall procedure is summarized in the following algorithm.

\begin{framed}
	\begin{textalgorithm}[Learning coarse-grained models from data]\label{algo:learning}
		\quad
		\begin{enumerate}[itemsep=0ex, topsep=0.5ex]
			\item Generate $m$ measurements $\{x_l\}_{l=1}^m$ of the aggregated state of the ABM.
			\item Estimate the drift and diffusion terms $\{b(x_l)\}_{l=1}^m$ and $ \{a(x_l)\}_{l=1}^m$ at the measurement points, e.g., via Monte Carlo simulations for short lag times $\tau$ using the Kramers--Moyal formulae \eqref{eq:kramersmoyal}.
			\item If applicable, normalize the data:
			\begin{enumerate}[itemsep=0ex, topsep=0ex]
				\item Reduce the training data by keeping only $d-1$ components of each measurement as well as its drift and diffusion estimates.
				\item Scale by the number of agents $N$, i.e., $c_i(t) = \frac{x_i(t)}{N}$.
			\end{enumerate}
			\item Choose a suitable set of basis functions $\{ \psi_i\}_{i=1}^n$ and compute the matrices $\Psi_X$ and $\mathrm{d}\Psi_X$.
			\item Minimize $\norm{\smash{\mathrm{d}\Psi_X - M \Psi_X}}_F $ and obtain a generator approximation $L=M^\top$ and identify the drift and diffusion terms using \eqref{eq:diff_identification}.
		\end{enumerate}
	\end{textalgorithm}
\end{framed}

\section{Numerical results}
\label{sec:numericalresults}

We will now apply Algorithm \ref{algo:learning} to three benchmark problems. First, we compare the numerical result with the theoretical SDE limit model \eqref{eq:SDE_limit} for the EVM in Section~\ref{exp:guiding} for varying numbers of agents $N$  and numbers of Monte Carlo samples $k$ for the pointwise drift and diffusion estimates as these are two crucial parameters for the quality of the numerically obtained model. In Section~\ref{sec:clustered}, we will then show that it can also be applied to the case where the network is not fully connected but consists of clusters connected by a few edges only. In Section~\ref{sec:PPM_results} we show for the PPM that it is also possible to obtain a reduced model for systems not based on interaction networks.

All results are compared using the \emph{root mean square error} (RMSE), which is defined by
\begin{equation*}
	\text{err} \coloneqq \left( \frac{1}{l} \sum_{i=1}^{l} (\widehat{y}_i - y_i)^2 \right)^{1/2},
\end{equation*}
where $y_i$ and $\widehat{y}_i$ denote the measured quantity and its prediction, respectively.

\subsection{Complete networks}
\label{sec:complete_networks}

Let us consider the EVM defined in Section~\ref{exp:guiding} and assume that the network is complete. 
The state space of this ABM is given by the $d-1$ dimensional simplex $\mathbb{X}_N$, with
\begin{equation*}
	\mathbb{X}_N \coloneqq \left\{x \in \mathbb{N}^d_0 \colon \sum_{i=1}^{d} x_i = N\right\}.
\end{equation*}
We consider now $d=3$ types and set the rate constants to
\begin{subequations}\label{eq:rate_constants_expl}
	\begin{align}
		\gamma_{12} = \gamma_{23} = \gamma_{31} &= 2 \label{eq:rate_constants_expl_imi1}, \\
		\gamma_{32} = \gamma_{21} = \gamma_{13} &= 1, \label{eq:rate_constants_expl_imi2} \\
		\gamma_{ij}' &= 0.01,
	\end{align}
\end{subequations}
for $i,j = 1,\dots,3$. Due to the conservation law, this is essentially a two-dimensional system. Thus, we eliminate one equation of the limit SDE \eqref{eq:SDE_limit} such that we can compare it with the data-driven SDE obtained by Algorithm \ref{algo:learning}. Additionally, after scaling the measurements by the number of agents, $N$, we obtain
\begin{equation}\label{eq:reduction}
	c_j(t) = 1 - \sum_{i \neq j} c_i(t).
\end{equation}
We will then evaluate the quality of the identified coarse-grained model.

Utilizing $c_3(t) = 1 - c_1(t) - c_2(t)$, we obtain the drift and diffusion terms
\begin{subequations}\label{eq:2dSystem}
	\begin{align}
		b & \colon \mathbb{X} \to \mathbb{R}^2, \label{eq:2D_drift_a}\\
		a & \colon \mathbb{X} \to \mathbb{R}^{2\times 2}, \label{eq:2D_diff_a}
	\end{align}
\end{subequations}
respectively. Note that $a(c) = a(c)^\top = (a_{ij}(c))$. Their derivation can be found in Appendix~\ref{adpx:reducedsystem}.

Following the arguments in Section~\ref{sec:learningSDE}, for a correct identification, we need a set of basis functions comprising monomials up to degree 3 as the highest order transition is of order $2$. For any given number of agents $N$, we can construct the first columns of the approximation $L_N$ of the generator $\mathcal{L}$ analytically via the coefficients of $b$ and $a$. E.g., for $N=10$ we obtain the matrix entry $l_{22}$ from the coefficient of $c_1$ in $b_1$, i.e.,  $l_{22} = \gamma_{31} - \gamma_{13} - \gamma_{12}' - \gamma_{13}' - \gamma_{31}'$, see \eqref{eq:red_drift} and \eqref{eq:red_diffusion} in Appendix~\ref{adpx:reducedsystem} for details. The first columns of $L_{10}$ are then given by
\begin{equation*}
	L_{10} =
	\kbordermatrix{
		         & 1 & c_1  & c_2    & c_1^2 & c_1c_2 & c_2^2 & \dots \\
		1        & 0 & 0.01 & 0.01   & 0.001 & 0     & 0.001 & \dots  \\
		c_1      & 0 & 0.97 & 0      & 0.321 & 0.009 & 0     & \dots  \\
		c_2      & 0 & 0    & -1.03  & 0     & 0.009 & 0.321 & \dots  \\
		c_1^2    & 0 & -1   & 0      & 1.64  & 0     & 0     & \dots  \\
		c_1c_2   & 0 & -2   & 2      & 0     & -0.36 & 0     & \dots  \\
		c_2^2    & 0 & 0    & 1      & 0     & 0     & -2.36 & \dots  \\
		c_1^3    & 0 & 0    & 0      & -2    & 0     & 0     & \dots  \\
		c_1^2c_2 & 0 & 0    & 0      & -4    & 1     & 0     & \dots  \\
		c_1c_2^2 & 0 & 0    & 0      & 0     & -1    & 4     & \dots  \\
		c_2^3    & 0 & 0    & 0      & 0     & 0     & 2     & \dots
	} \in \R^{10 \times 10}.
\end{equation*}

We will compare the numerical results to the corresponding columns of $L_N$ and the drift and diffusion terms \eqref{eq:2D_drift_a} and \eqref{eq:2D_diff_a}, respectively. The identified system has the following structure:
\begin{subequations}\label{eq:DDSystem}
	\begin{align}
		b_i(c) &\coloneqq \beta^i_5 \ts c_1^2 + \beta^i_4 \ts c_2^2 + \beta^i_3 \ts c_1c_2 + \beta^i_2 \ts c_1 + \beta^i_1 \ts c_2 + \beta^i_0 \label{eq:reduced_drift_b}, \\
		a_{ij}(c) &\coloneqq \kappa^{ij}_5 \ts c_1^2 + \kappa^{ij}_4 \ts c_2^2 + \kappa^{ij}_3 \ts c_1c_2 + \kappa^{ij}_2 \ts c_1 + \kappa^{ij}_1 \ts c_2 + \kappa^{ij}_0, \label{eq:reduced_diff_b}
	\end{align}
\end{subequations}
where the coefficients are given by the expressions derived for \eqref{eq:2D_drift_a} and \eqref{eq:2D_diff_a}, see Appendix~\ref{adpx:reducedsystem}. The coefficients $\beta_h^i$ of \eqref{eq:reduced_drift_b} can immediately be obtained from the second and third column of $L_N$. The coefficients $\kappa_{h}^{ij}$ of \eqref{eq:reduced_diff_b} are extracted from the columns four to six by using \eqref{eq:diff_identification}. E.g., for $a_{12}(c)$ we obtain
\begin{alignat*}{2}
	b_1(c) &= (\mathcal{L}\psi_2)(c) &&=  - c_1^2 - 2 \ts c_1c_2 + 0.97 \ts c_1 + 0.01,\\
	b_2(c) &= (\mathcal{L}\psi_3)(c) &&=  c_2^2 +2 \ts c_1c_2 - 1.03 \ts c_2 + 0.01,\\
	a_{12}(c) &= (\mathcal{L}\psi_5)(c) - b_1(c) \ts c_2 - b_2(c) \ts c_1 &&= - 0.3 \ts c_1c_2 - 0.001 \ts c_1 - 0.001 \ts c_2.
\end{alignat*}

Comparing the coefficients of the SDE limit model with its corresponding parts in the data-driven system, we can (under certain conditions) recover the rate constants of the underlying Markov jump process. For the considered example we compare the coefficients of \eqref{eq:2dSystem} with \eqref{eq:DDSystem}. We set up a system of linear equations $A \gamma = v$ for a suitable matrix $A$, where $\gamma$ and $v$ are given by
\begin{align*}
	\gamma &= [\gamma_{12}, \ts \gamma_{13}, \dots, \gamma_{32}]^\top, \\
	v &= [\beta_5^1, \dots, \beta_0^2, \ts\kappa_5^{11}, \dots, \kappa_0^{22}]^\top.
\end{align*}
Note that for this example with the rate constants chosen in \eqref{eq:rate_constants_expl} the system cannot be solved exactly in general since the model is symmetric in the sense that imitation is possible in both ways (i.e., $\gamma_{ij} \neq 0$ for all $i\neq j$). Thus, we only find values for $\gamma_{ij}$ and $\gamma_{ji}$ satisfying the differences appearing in \eqref{eq:2D_drift_a} and \eqref{eq:2D_diff_a}, see Appendix~\ref{adpx:reducedsystem} for details. However, this has only an influence on the reconstruction of the underlying Markov jump process but not on the coarse-grained model.

\subsubsection*{Evaluations}

For both the number of agents $N$ and the number of Monte Carlo samples $k$, we set a maximum of 5000. Since the state space $\mathbb{X}_N$ is discrete and $N$ constant, the amount of distinct points is finite and depends on $N$ and $d$; more precisely for a $d$-dimensional regular discrete simplex with $N+1$ points on each edge, the number of points is given by $\binom{N+d}{d}$ for $d \le N$ \cite{Costello71}. In our example, we have a two-dimensional simplex and thus $\binom{N+2}{2}$ points. The number of uniformly chosen measurements is given in Table~\ref{tab:my-table} for different $N$. We then estimate the drift and diffusion term for each point via~\eqref{eq:kramersmoyal} for $k$ short simulations of the Markov jump process with a lag time of $\tau=0.01$ resulting in a total of $m \cdot k$ training data points.

\begin{figure}[!t]
	\centering
	\begin{subfigure}[t]{0.49\textwidth}
		\centering
		\caption{Drift term}
		\includegraphics[scale=0.9]{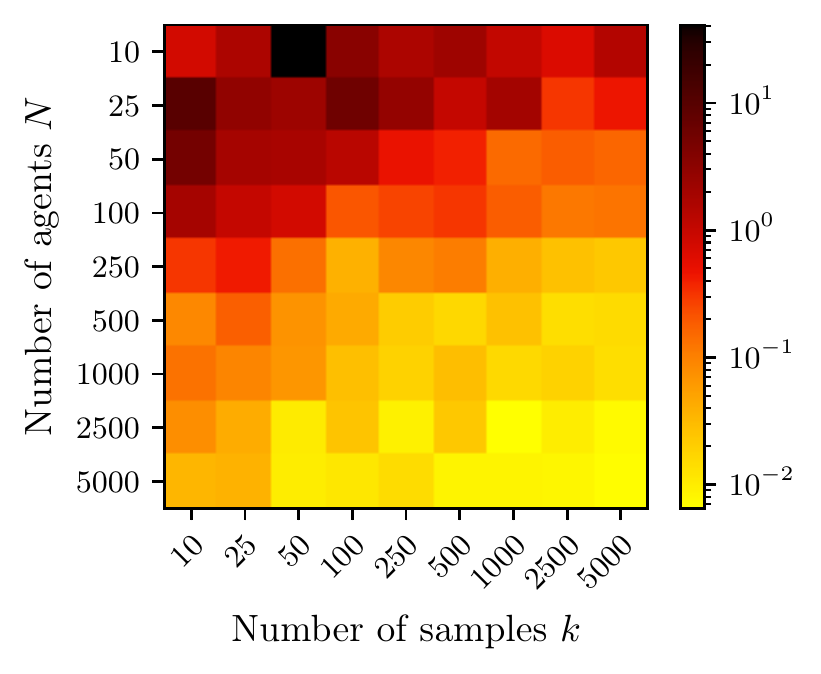}
	\end{subfigure}
	\hfill
	\begin{subfigure}[t]{0.49\textwidth}
		\centering
		\caption{Diffusion term}
		\includegraphics[scale=0.9]{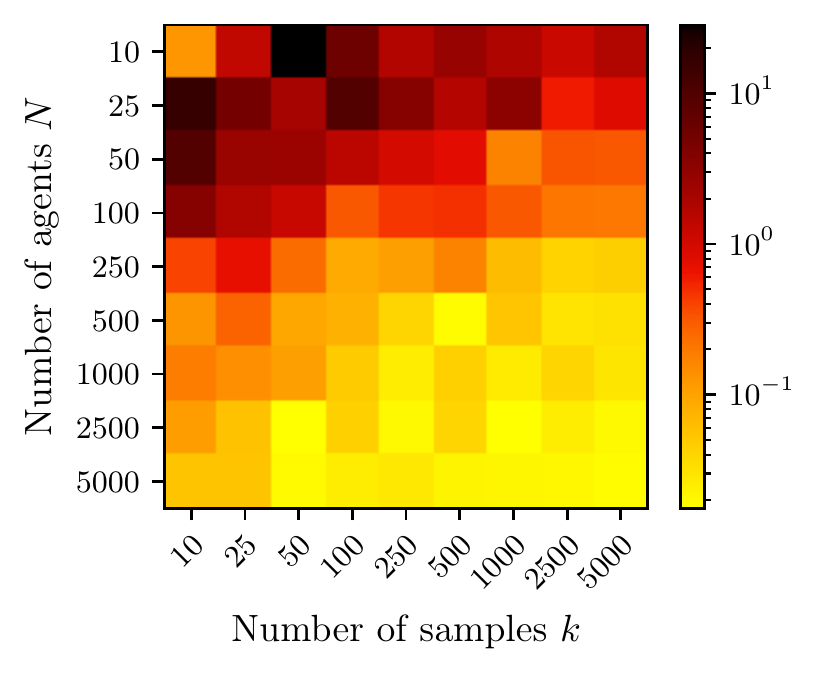}
	\end{subfigure}
	\caption{Approximation error defined as the RMSE of the coefficients of (a) the drift and (b) diffusion estimates for the EVM in Section~\ref{exp:guiding} compared to the exact SDE limit model \eqref{eq:SDE_limit} depending on the number of agents $N$ and number of Monte Carlo samples $k$ for the estimation via Kramers--Moyal formulae. The brighter the color, the smaller the error and the better the identification of the reduced system. For increasing $N$ and $k$ the approximation error decreases.}
	\label{fig:errCoeffThree}
\end{figure}

Figure~\ref{fig:errCoeffThree} shows the approximation error of the numerically obtained coefficients and their theoretical counterparts  appearing in \eqref{eq:2D_drift_a} and \eqref{eq:2D_diff_a} depending on the number of agents and the number of Monte Carlo samples. For both parameters, the RMSE decreases by several orders of magnitude as $N$ and $k$ increase. Note that the number of agents $N$ has a significantly larger influence than the number of samples $k$. Especially for small $N$, e.g., $N=10$, we observe that higher values of $k$ do not improve the results. This is consistent with the literature as the SDE model \eqref{eq:SDE_limit} approximates the Markov jump process for large~$N$.

\begin{figure}[!t]
	\centering
	\begin{subfigure}[t]{0.45\textwidth}
		\centering
		\caption{}
		\includegraphics[scale=0.9]{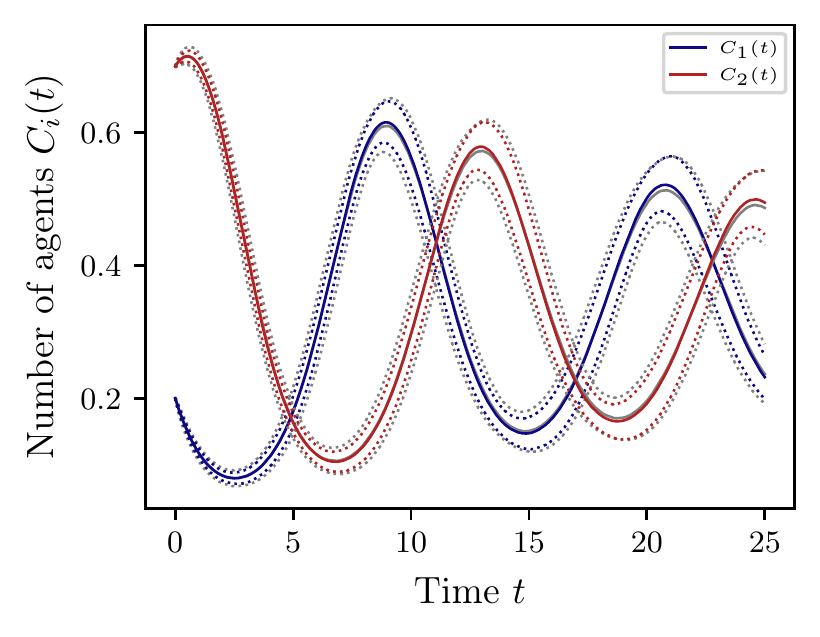}
	\end{subfigure}
	\hfil
	\begin{subfigure}[t]{0.45\textwidth}
		\centering
		\caption{}
		\includegraphics[scale=0.9]{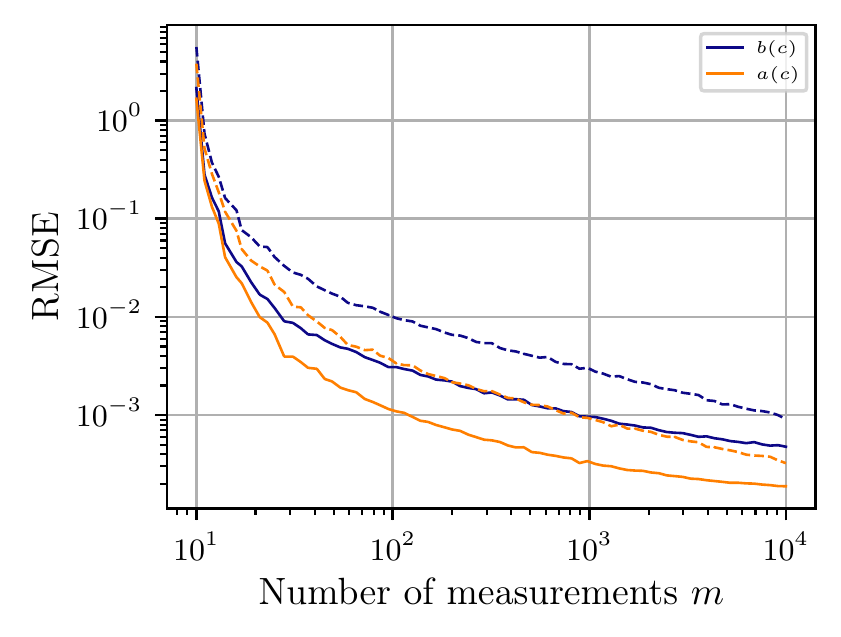}
	\end{subfigure}
	\caption{(a) Expectation (solid) and standard deviation (dashed) of the SDE limit model $C_i(t)$ and its data-driven approximation (gray) estimated from $10^3$ Monte Carlo simulations for the dynamics of the EVM of Section~\ref{exp:guiding} for $N=5000$ agents and initial state $c(0) = [0.2, 0.7, 0.1]^\top \in \mathbb{X}$. The relative number of agents of type $S_3$ can be reconstructed using \eqref{eq:reduction} and is therefore not displayed. The approximate moments (gray solid and dashed lines) agree with the SDE limit model. (b) Approximation and evaluation error of the drift and diffusion estimates for the EVM in Section~\ref{exp:guiding} compared to the exact SDE limit model \eqref{eq:SDE_limit} depending on the number of measurements $m$ for fixed $k_1=10$ (dashed), $k_2=100$ (solid) and $N=5000$ agents. The error is averaged over 100 simulations. Clearly, for higher amounts of training data a smaller error can be expected. This holds for both parameters $m$ and $k$.}
	\label{fig:longprediction}
\end{figure}

As it is not only important to identify the coefficients of an SDE limit model, we also compare how well the reduced model approximates the dynamics of the ABM, e.g., to make predictions about the number of agents of a specific type. Figure~\ref{fig:longprediction}~(a) shows a comparison for a long-time realization in terms of expectation (solid line) and standard deviation (dashed line) for the data-driven model and its theoretical equivalent estimated from $10^3$ Monte Carlo samples. Both first- and second-order moments are almost indistinguishable from the theoretical SDE limit model. Considering the numerical effort that renders the simulation required for Figure~\ref{fig:longprediction}~(a) infeasible in many cases, the estimated coarse-grained model yields valuable results. Additionally, it is obtained in a fraction of the time it takes to simulate the original ABM.

Figure~\ref{fig:longprediction}~(b) shows the dependency of the RMSE on the number of measurements $m$ for two fixed $k_i$, namely $k_1=10$ (dashed line) and $k_2=100$ (solid line). The error is averaged over 100 simulations for $5000$ agents. We observe that for greater $m$ the error, as expected, decreases by several orders of magnitude, independently of $k$. However, the impact of increasing $m$ is larger than the one of increasing $k$. For small values $m \cdot k_i$, the error is smaller for $k_1 = 10$ (dashed line) than for $k_2=100$ because the measurements cover the state space more densely: For $k_1 = 10$, for example, we have $m=10$ measurements while for $k_2 = 100$ we only have $m=1$ measurement. Thus, there are two tuning parameters for the amount of training data to be used.

\begin{table}[!h]
	\centering
	\caption{Measurement set sizes for given number of agents $N$. For each measurement, we estimate the drift and diffusion term with $k$ short Monte Carlo simulations for a lag time $\tau = 0.01$ resulting in a total training data set size of $m\cdot k$.}
	\label{tab:my-table}
	\begin{tabular}{rr}
		\toprule
		Number of agents $N$ & Measurements $m$      \\ \midrule
		10                            & 7          	 \\
 		25                            & 35        	 \\
		50                            & 133       	 \\
		100                           & 515          \\
		250                           & 3163      	 \\
		$N\ge500$                     & 10000     	 \\   \bottomrule
	\end{tabular}
\end{table}

\subsection{Clustered networks}
\label{sec:clustered}

Let us now consider the case where the network consists of $Q$ (not necessarily equally-sized) clusters. Within a cluster each agent is connected to all other agents, i.e., each cluster $q$ is a complete sub-graph of size $N_q$. Two agents of different clusters are connected with probability $p$. If $p$ is sufficiently small, then the clusters are connected only by a few edges and the corresponding sub-matrix of the adjacency matrix is sparse. As before, each agent is influenced by its neighbors. However, due to the non-completeness of the network, the resulting transition propensities depend on the size of the individual neighborhood; therefore, they might differ among agents. Here, we do not model the population state of the ABM as described in Section~\ref{sec:Modelingagentbasedsystems} since the overall aggregation leads to errors in this case. Instead, we augment the population state by subpopulations, i.e., an aggregation by cluster. We will use these to learn a coarse-grained model of the agent dynamics.

\subsubsection*{An SDE limit model for clustered networks}

We can set up a limit model that describes the relative frequencies of each type per cluster. As mentioned before, this limit model contains an approximation error that is due to the aggregation of types in each cluster. However, under certain conditions (e.g., uniformly drawn connecting edges) the model yields a good approximation.

We extend \eqref{eq:SDE_limit} such that it describes the temporal evolution of the relative frequencies for a network that consists of $Q$ clusters. Assume that the connecting edges are drawn uniformly with probability $p$. Let $N$ be the number of agents in cluster $q = 1, \dots, Q$. For simplicity we assume that all clusters are equally sized. We augment the system state such that it has the relative frequencies of each type per cluster, i.e.,
\begin{equation*}
	C(t) = \left[c_1(t)^\top, \dots, c_Q(t)^\top  \right]^\top \in \mathbb{R}^{d \ts Q}.
\end{equation*}
Let $\widetilde{\alpha}_{q, k}$ be the rescaled propensity function for transition $k$ in cluster $q$ and $\nu_{q,k} \in \mathbb{R}^{d \ts Q}$ its corresponding net change vector. We obtain
\begin{equation}\label{eq:cluster_SDE}
	\mathrm{d}C (t) = \sum_{q=1}^{Q} \left[\ts \sum_{k=1}^{K_q} \widetilde{\alpha}_{q,k}(C(t)) \ts \nu_{q,k} \ts \mathrm{d}t + \sum_{k=1}^{K_q} \frac{1}{\sqrt{N}} \sqrt{\widetilde{\alpha}_{q,k}(C(t))} \ts \mathrm{d}W_{q,k}(t) \ts \nu_{q,k}\right].
\end{equation}
Note that equation \eqref{eq:cluster_SDE} can be rewritten so that clusters can also have different sizes. Given a cluster $q$ the diffusion term corresponding to transitions within the cluster scales with $1/\sqrt{N_q}$ while transitions induced by another cluster $q'$ scales with $1/\sqrt{N_{q'}}$.

\begin{expl}[SDE limit model for two clusters]
	Let us consider a network consisting of $Q=2$ clusters each having $N_1$ and $N_2$ agents, and let $p$ be the probability for an edge connecting two agents of cluster $Q_1$ and $Q_2$. We define the connection strength of cluster $Q_1$ and $Q_2$ as the ratio between the number of edges $E$ connecting both clusters and the total number of possible connecting edges $E_\text{max} =N_1N_2$. The expected connection strength is given by $p$ since
	\begin{equation*}
		\mathbb{E}\left[\frac{E}{E_\text{max}}\right] = \frac{\mathbb{E}[E]}{E_\text{max}} = \frac{p \ts N_1N_2}{N_1N_2} = p.
	\end{equation*}	
	As in Section~\ref{exp:guiding}, we consider imitation and exploration. The latter is independent of the considered network, while the former is either induced from the inside or outside. If the transition is caused from the inside, we call it \emph{intra-cluster} transition and \emph{inter-cluster} transition if it is caused from the outside. Intra-cluster transitions are denoted by $R_{ij}$ and $R_{ij}'$. Imitation as an inter-cluster transition rule is given by
	\begin{equation*}
		R_{qq', ij} \colon S_{q,i} + S_{q',j} \to S_{q,j} + S_{q',j}.
	\end{equation*}
	For the intra-cluster transitions the propensity functions are given by
	\begin{equation*}
		\alpha_{ij} = \frac{1}{N_q+p\ts N_{q'}}\gamma_{q,ij} \ts x_{q,i} \ts x_{q,j},
	\end{equation*}
	while for the inter-cluster transition they are given by
	\begin{equation*}
		\alpha_{qq',ij} = p \ts\frac{1}{N_q+p\ts N_{q'}}\beta_{q,ij} \ts x_{q,i} \ts x_{q',j}
	\end{equation*}
	as each agent has $N_q+p\ts N_{q'}$ possible partners for interaction.
	
	For simplicity, we assume that both clusters are of the same size. For the corresponding net change vector, it holds that $\nu_{qq',ij} = \nu_{q,ij}$ as the inter-cluster transitions $R_{qq', ij}$ only influences state $c_q(t)$ and not $c_{q'}(t)$. For $C(t)=\left[c_1(t)^\top,c_2(t)^\top\right]^\top \in \mathbb{R}^{2d}$, the SDE solution is given by
  \begin{footnotesize}
	\begin{subequations} \label{eq:two_cluster_SDE}
		\begin{align}
			\mathrm{d}c_{q,i}(t) = \Bigg[ & \sum_{i \neq j} \frac{1}{(p+1)} (\gamma_{q,ji}-\gamma_{q,ij}) \ts c_{q,i}(t) \ts c_{q,j}(t) \\
			+ &\sum_{i\neq j}\Big[ \gamma'_{q,ji} \ts c_{q,j}(t)-\gamma'_{q,ij} \ts c_{q,i}(t) \Big] \\
			+ &\sum_{i\neq j} \frac{p}{(p+1)} \Big[ \beta_{q,ji} \ts c_{q,j}(t) \ts c_{q',i}(t) - \beta_{q,ij} \ts  c_{q,i}(t) \ts c_{q',j}(t) \Big] \Bigg] \mathrm{d} t \\
			+ \ts \frac{1}{\sqrt{N}} \Bigg[ &\sum_{i \neq j} \sqrt{\frac{1}{(p+1)} \gamma_{q,ji} \ts c_{q,i}(t) \ts c_{q,j}(t)} \ts \mathrm{d}W_{q,ji}^{\mathrm{im}}(t) - \sqrt{\frac{1}{(p+1)} \gamma_{q,ij} \ts c_{q,i}(t) \ts c_{q,j}(t)} \ts \mathrm{d}W_{q,ij}^{\mathrm{im}}(t) \\
			+ &\sum_{i \neq j} \sqrt{\gamma'_{q,ji} \ts c_{q,j}(t)} \ts \mathrm{d}W_{q,ji}^{\mathrm{ex}}(t) - \sqrt{\gamma'_{q,ij} \ts c_{q,i}(t) } \ts \mathrm{d}W_{q,ij}^{\mathrm{ex}}(t) \\
			+ & \sum_{i\neq j} \sqrt{\frac{p}{(p+1)} \beta_{q,ji} \ts c_{q,j}(t) \ts c_{q',i}(t)} \ts \mathrm{d}W_{q,ji}^{\mathrm{int}}(t) - \sqrt{\frac{p}{(p+1)} \beta_{q,ij} \ts c_{q,i}(t) \ts c_{q',j}(t)} \ts \mathrm{d}W_{q,ij}^{\mathrm{int}}(t) \Bigg].
		\end{align}
	\end{subequations}
    \end{footnotesize}
	The addends (a), (b), (d), and (e) correspond to intra-cluster transitions, while (c) and (f) correspond to inter-cluster transitions. We will drop the index $q$ whenever it is clear from the context.  \exampleSymbol
\end{expl}

\subsubsection*{Evaluations}

\begin{figure}[!t]
	\centering
	\begin{subfigure}[t]{0.45\textwidth}
		\centering
		\caption{}
		\includegraphics[scale=0.9]{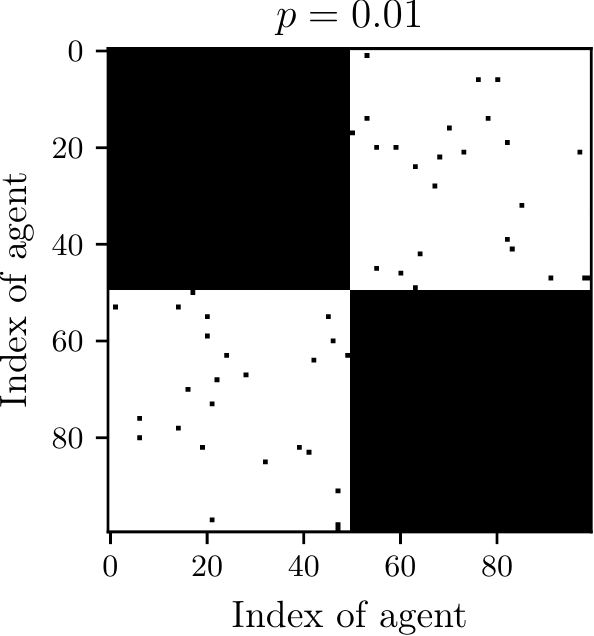}
	\end{subfigure}
	\hfill
	\begin{subfigure}[t]{0.45\textwidth}
		\centering
		\caption{}
		\includegraphics[scale=0.9]{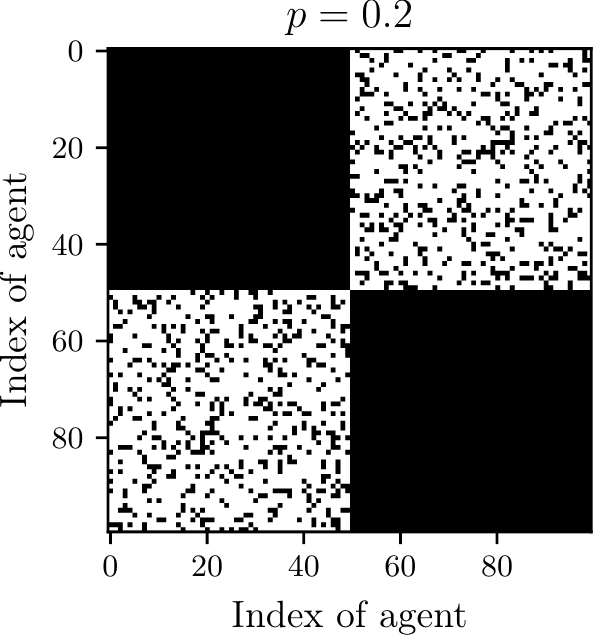}
	\end{subfigure}
	\begin{subfigure}[t]{0.45\textwidth}
		\centering
		\caption{}
		\includegraphics[scale=0.9]{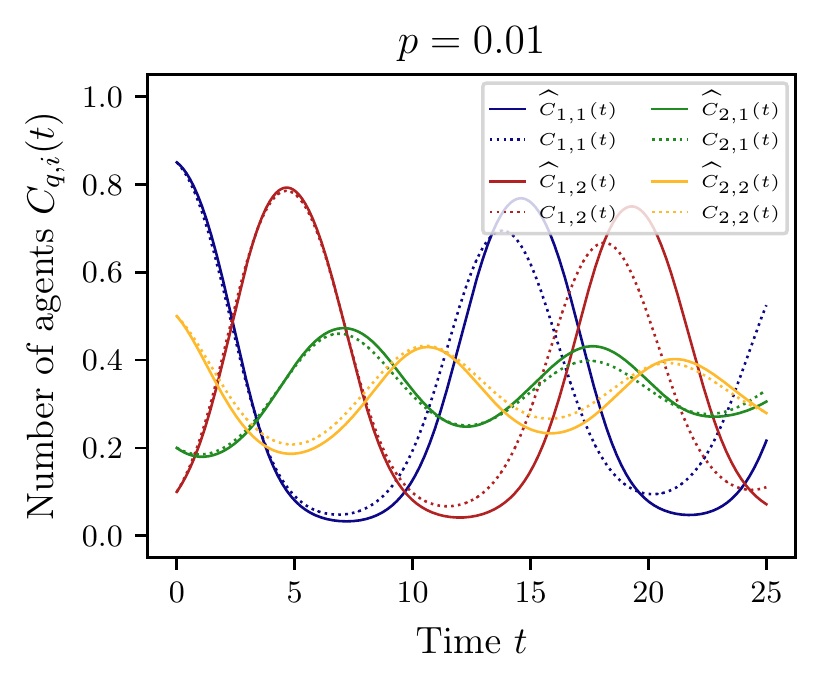}
	\end{subfigure}
	\hfill
	\begin{subfigure}[t]{0.45\textwidth}
		\centering
		\caption{}
		\includegraphics[scale=0.9]{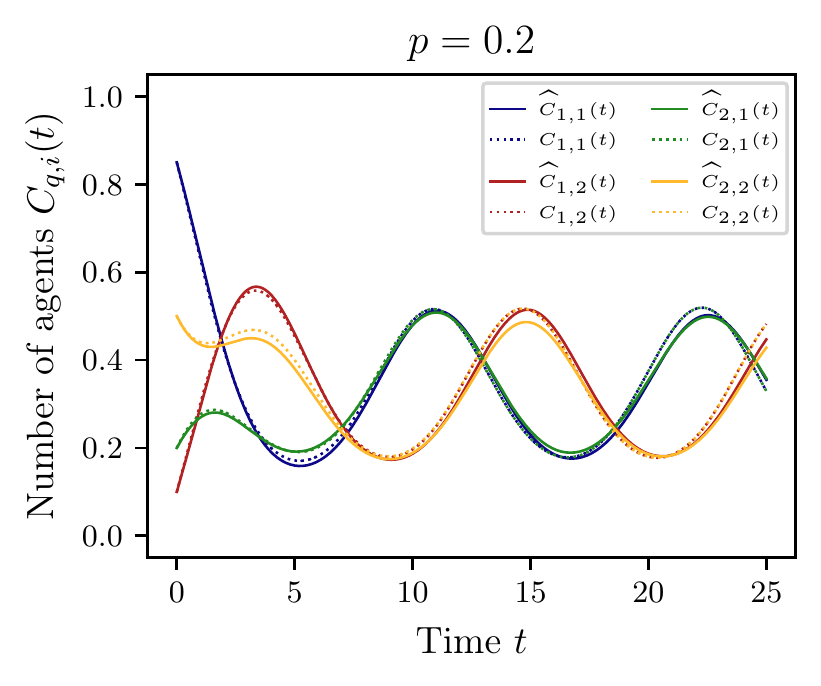}
	\end{subfigure}
	\caption{(a) \& (b) Adjacency matrices of the networks where black represents 1 (existing edge) and white 0 (no edge). (c) \& (d) First-order moment of the data-driven coarse-grained model (solid) and the limit SDE (dotted) \eqref{eq:two_cluster_SDE} for two clusters with $N=50$ agents, $\gamma_{12} = \gamma_{23} = \gamma_{31} = 2$, $\gamma_{13} = \gamma_{21} = \gamma_{32} = 1$, $\gamma_{ij}' = 0$ for all $i,j = 1, \dots, 3$ and $c(0) = [0.85, 0.1, 0.05, 0.2, 0.5, 0.3]^\top$. The data-driven model is estimated using $k=1000$ realizations of $m=1000$ measurements for lag time $\tau = 0.01$.}
	\label{fig:clustered_solution}
\end{figure}

We now simulate the EVM in discrete time with step size $t_\text{step} = 0.01$, see Appendix~\ref{apdx:DTEVM} for the pseudocode. While it can be applied to arbitrary networks, we restrict ourselves to highly clustered networks as depicted in Figure~\ref{fig:clustered_solution}~(a) and (b). We create $k=1000$ realizations for each of the $m=1000$ uniformly drawn initial states of the ABM for a lag time of $\tau = 0.01$. The network consists of two equally sized clusters, each containing $N=50$ agents. We assume $\gamma_{q,ij} = \gamma_{q',ij}$, $\gamma'_{q,ij} = \gamma'_{q',ij}$, and $\beta_{q,ij} = \beta_{q',ij} = \gamma_{q,ij}$ for all $i, j$. The rate constants for imitative transitions are given by \eqref{eq:rate_constants_expl_imi1} and \eqref{eq:rate_constants_expl_imi2}. For exploratory transitions we set $\gamma_{ij}' = 0$ for all $i,j$.

We compare the data-driven model and the model defined in \eqref{eq:two_cluster_SDE} for two networks with different connection strengths. The adjacency matrices of both networks are shown in Figure~\ref{fig:clustered_solution}~(a) and (b). The first network has a connection strength of $p=0.01$ while the second has a 20-times larger connectivity, i.e., $p=0.2$. The first network is a subgraph of the second. We apply Algorithm \ref{algo:learning} to the cluster-based aggregate states of the agent dynamics for each network to obtain the data-driven coarse-grained model. Figure~\ref{fig:clustered_solution}~(c) and (d) show the prediction of the temporal evolution of the first moments for each type per cluster. Note that the colors are different from Figure~\ref{fig:guiding}. Both realizations start from the same initial value. The difference in their temporal evolution results directly from the network structure. As described in Section~\ref{sec:complete_networks} for complete networks, the results improve for larger values of $N$, $m$, and $k$. We can also observe in Figure~\ref{fig:clustered_solution}~(d) that for a higher connectivity, i.e., larger $p$, both clusters synchronize so that the relative numbers of agents per type are identical in each cluster.

\begin{figure}[!t]
	\centering
	\includegraphics[scale=0.9]{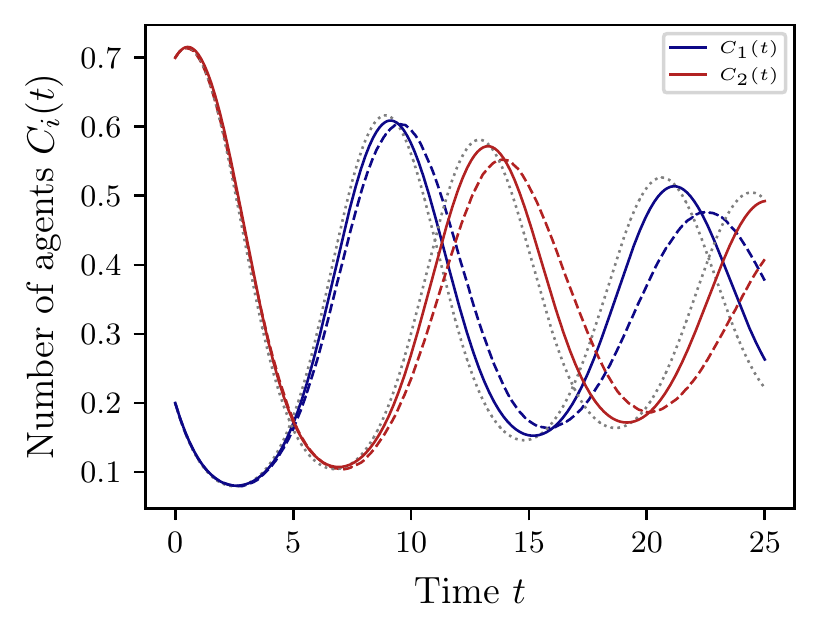}
	\caption{Expectation of the data-driven reduced model (solid) compared the EVM (dashed) on a random network with average degree of approximately $50$, estimated from $10^3$ Monte Carlo simulations for $N=500$ agents and initial state $c(0) = [0.2, 0.7, 0.1]^\top \in \mathbb{X}$. The deterministic part of the SDE limit model \eqref{eq:SDE_limit} is indicated in gray (dotted). The data-driven model is estimated using $m=k=1000$ measurements and realizations for the lag time $\tau = 0.01$.}
	\label{fig:incompletegedmd}
\end{figure}

\begin{rem}
	Consider a random network of $N=500$ nodes where two agents are connected with a probability of 10~\%. The resulting network is sparsely connected and exhibits an approximate average degree of 50. Figure~\ref{fig:incompletegedmd} shows the expectation of the data-driven model (solid) compared to the EVM (dashed) for this random network estimated from $10^3$ Monte Carlo simulations. For short times $t$, the data-driven model agrees with the ABM. However, for larger time $t$ the prediction deteriorates mainly due to the sparsity of the network. Note that the absence of a ground truth model for the EVM on this sparse network complicates the analysis. Compared to the expectation obtained via the SDE limit model~\eqref{eq:SDE_limit} (indicated in gray, dotted) the data-driven model yields a better approximation.
\end{rem}

\subsection{Predator-prey model}
\label{sec:PPM_results}

\begin{figure}[!t]
	\centering
	\begin{subfigure}[t]{0.32\textwidth}
		\centering
		\caption{}
		\includegraphics[width=\textwidth]{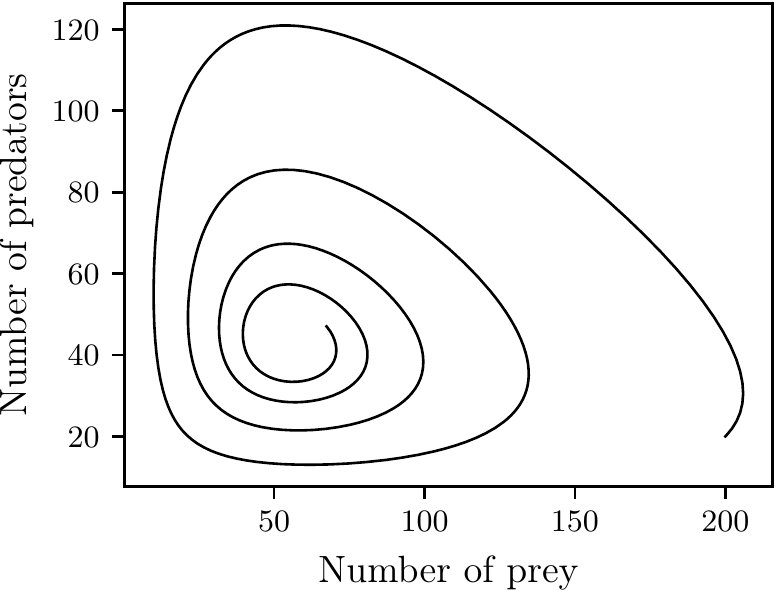}
	\end{subfigure}
	\hfill
	\begin{subfigure}[t]{0.32\textwidth}
		\centering
		\caption{}
		\includegraphics[width=\textwidth]{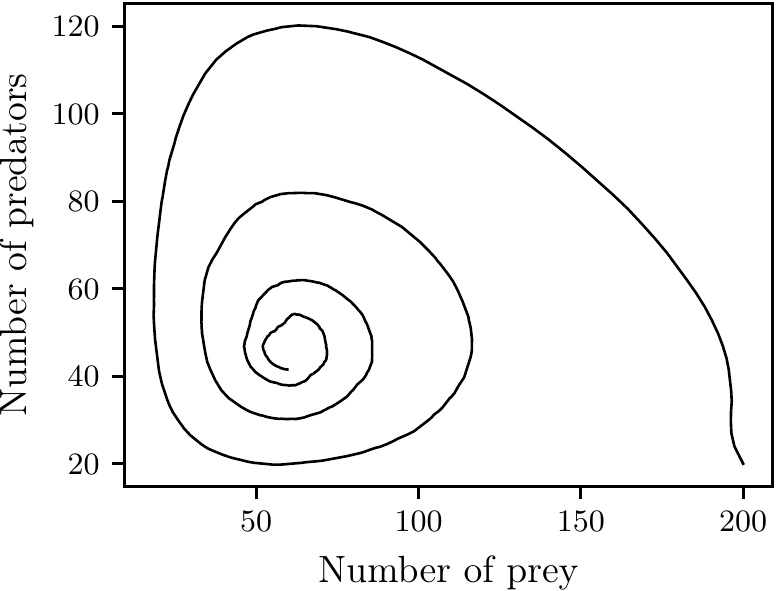}
	\end{subfigure}
	\hfill
	\begin{subfigure}[t]{0.32\textwidth}
		\centering
		\caption{}
		\includegraphics[width=\textwidth]{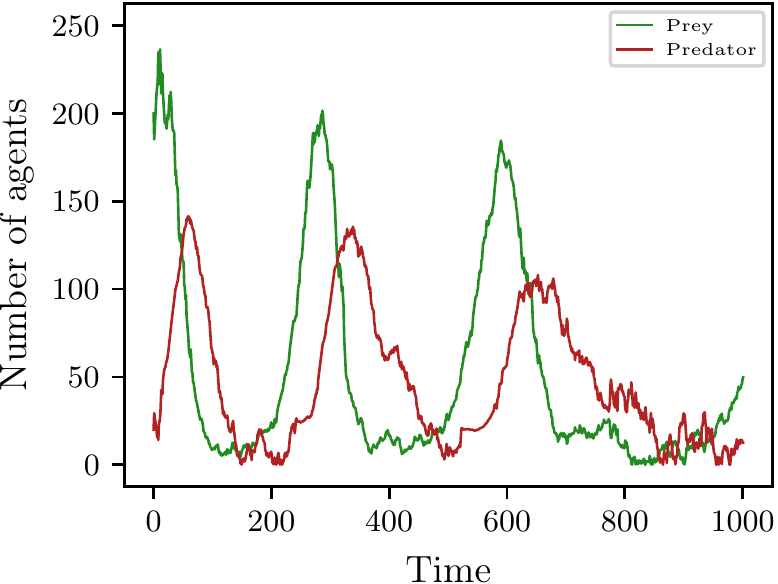}
	\end{subfigure}
	\caption{Phase portraits of first-order moment of (a) the reduced SDE model and (b) the PPM estimated from 958 Monte Carlo simulations. (c) Realization of the reduced SDE solution learned from $m=k=1000$ measurements and samples for the PPM with parameters given in Table \ref{tab:PPM_parameters}.}
	\label{fig:PPM_reduced}
\end{figure}

Let us now consider the PPM introduced in Section~\ref{sec:PPM}. The parameters we use are listed in Table~\ref{tab:PPM_parameters}. We learn a data-driven model from $m=k=1000$ measurements and samples. The lag time for estimating drift and diffusion is set to $\tau = 1$. Although the defined PPM has a spatial component, i.e., relatively slow movement of the agents with respect to the dimension of the space and search radius $v$ of the predators, we use the classic Lotka--Volterra differential equations as a starting point for the set of basis functions. The set consists of monomials up to degree 3 so that we can identify the coefficients of the drift and diffusion terms. Figure~\ref{fig:PPM_reduced}~(a) and (b) show the phase portrait of the first-order moment of the reduced SDE model and the PPM averaged over 958 realizations. In 42 out of 1000 realizations the predators died out before the prey so that the size of the prey population grows exponentially. The results show that the reduced model is able to approximate the qualitative dynamical behavior of the PPM. Figure~\ref{fig:PPM_reduced}~(c) shows a realization of the reduced SDE model.

\begin{table}[!h]
	\centering
	\caption{Parameters used during the simulation of the PPM}
	\label{tab:PPM_parameters}
	\begin{tabular}{ll}
		\toprule
		Parameter                            					& Value             \\ \midrule
		Space height $\times$ width                  			& 100 $\times$ 100  \\
		Variance $h$	                   						& 1  				\\
		Reproduction probability prey $p_\text{rep}$     	  	& 0.03              \\
		Reproduction probability predator	$p_\text{rep}'$ 	& 0.5				\\
		Probability of death $p_\text{death}$					& 0.02				\\ 
		Radius of vision $v$					 				& 3					\\ \bottomrule
	\end{tabular}
\end{table}

\section{Conclusion}
\label{sec:conclusion}

In this work, we showed how the Koopman generator can be used to obtain coarse-grained stochastic models from aggregate state data of agent-based dynamics. We demonstrated the procedure for two different ABMs, namely a voter model and a predator-prey model. The ABM codes used for generating the results presented in this paper can be found at \url{https://github.com/Henningston/ABMs}.

In the first case we considered complete and clustered interaction networks of homogeneous agents such that each agent can interact at any time with all other agents (or within their cluster, respectively). We showed that under certain conditions the reduced models agree with their respective SDE limit models. In both considered cases, we showed that the data-driven reduced models are suitable for predictions. The results of Section~\ref{sec:complete_networks} showed that when considering incomplete, clustered networks, aggregation of state variables led to an approximation error in the population state model. As a consequence, the data-driven model and its SDE approximation agreed only for short time intervals, see Figure~\ref{fig:clustered_solution}~(a). It also showed that the number of agents per cluster needs to be large enough or, alternatively, the connectivity between them high enough for the data-driven coarse-grained model and the SDE model \eqref{eq:two_cluster_SDE} to agree, see Figure~\ref{fig:clustered_solution}~(b). First experiments showed that for networks with an arbitrary structure the prediction horizon can be shorter which implies that, if the state of an ABM depends strongly on the spatial structure, e.g., formation of clusters, coexistence or spatial heterogeneity, this needs to be taken into account, see Figure~\ref{fig:incompletegedmd}.

For the second model -- the predator-prey system -- we showed in Section~\ref{sec:PPM_results} that it is also possible to identify a reduced model for an ABM that is not bound to interaction networks and whose time step is comparably large (i.e., not close to zero as in the first case). The reduced model is able to capture the qualitative behavior.

Our approach is limited to ABMs where it is believed that the aggregated dynamics can be meaningfully represented by ODEs or SDEs. However, this approach might fail if spatial interaction or interaction with the space itself have a strong influence on the behavior of the agents and therefore the outcome of the model.

In general, our approach relies on the assumption that all types of agents are available in sufficient numbers. If the number of agents (more generally speaking the size of the system) is large enough, it is known that the SDE accurately approximates the chemical master equation \cite{grima2011accurate}. However, there exist cases where the SDE fails to capture the behavior of a discrete ABM, more precisely noise-induced metastability. This is the case when bi- or multi-stability stems from the discreteness of the system (that is, if the size of the system is not large enough) \cite{duncan2015noise,hanggi1984bistable}; see also Figure~\ref{fig:guiding}~(b) and Figure~\ref{fig:longprediction}~(a) for systems with small and large numbers of agents, respectively.

Additionally, the approach relies on accurate, \emph{pointwise} estimates of the drift and diffusion terms. Inaccurate, insufficient estimates lead to nonsparse solutions of the generator approximation. Additional techniques like \emph{iterative hard thresholding} or \emph{denoising} might be applied to improve the results, see \cite{Klus2020} and references therein.

Our approach to obtain data-driven coarse-grained models from agent-based dynamics opens up new possibilities for further analysis and has the potential to reduce the numerical effort when investigating ABMs. In addition to parameter optimization or sensitivity analysis, which are often infeasible due to the complexity of the ABM, the reduced model can also be used to find control schemes to steer the system to a desired state. More precisely, the reduced model can be used to find, e.g., harvesting schedules for systems like the predator-prey models or to develop strategies to persuade agents to change their opinion (e.g. electoral or commercial campaigns, or use of green technology). Future research will address the control of ABMs using data-driven reduced models.

\section*{Acknowledgements}

This research has been funded by Germany's Excellence Strategy (MATH\texttt{+}: The Berlin Mathematics Research Center, EXC-2046/1, project ID: 390685689) and through Deutsche Forschungsgemeinschaft (DFG, German Research Foundation) through grant CRC 1114 (Scaling Cascades in Complex Systems, project ID: 235221301). We acknowledge support by  the Open Access Publication Fund of the Freie Universität Berlin. We thank Stefanie Winkelmann for helpful discussions.

\bibliographystyle{plain}
\bibliography{bibliography}
	
\appendix
\newpage
\section*{Appendix}

\section{Reduced two-dimensional system}
\label{adpx:reducedsystem}

Given the extended voter model with $d=3$ types, we choose the rate constants
\begin{align*}
	\gamma_{12} &= \gamma_{23} = \gamma_{31} = 2, \\
	\gamma_{32} &= \gamma_{21} = \gamma_{13} = 1, \\
	\text{and} \quad\gamma_{ij}' &= 0.01 \quad \text{for} \quad i,j = 1,\dots,3.
\end{align*}
Due to the conservation law, this is essentially a two-dimensional system. Utilizing $c_3(t) = 1 - c_1(t) - c_2(t)$, we obtain the (reduced) drift term $b\colon \mathbb{X} \to \mathbb{R}^2$ given by
\begin{subequations}\label{eq:red_drift}
	\begin{align}
		\begin{split}
			b_1(c) &= (\gamma_{13} - \gamma_{31}) \ts c_1^2 + (\gamma_{21} - \gamma_{12} + \gamma_{13} - \gamma_{31}) \ts c_1c_2 \\
			&+ (\gamma_{31} - \gamma_{13} - \gamma_{12}' - \gamma_{13}' - \gamma_{31}') \ts c_1 + (\gamma_{21}' - \gamma_{31}') \ts c_2 + \gamma_{31}',
		\end{split}
		\\
		\begin{split}
			b_2(c) &=  (\gamma_{23} - \gamma_{32}) \ts c_2^2+ (\gamma_{12} - \gamma_{21} + \gamma_{23} - \gamma_{32}) \ts c_1c_2\\
			&+ (\gamma_{32} - \gamma_{23} - \gamma_{21}' - \gamma_{23}' - \gamma_{32}') \ts c_2 + (\gamma_{12}' - \gamma_{32}') \ts c_1 + \gamma_{32}'.
		\end{split}
	\end{align}
\end{subequations}
The (reduced) diffusion term $a\colon \mathbb{X} \to \mathbb{R}^{2\times 2}$, $a(c) = a(c)^\top = (a_{ij}(c))$ is given by
\begin{subequations}\label{eq:red_diffusion}
	\begin{align}
		\begin{split}
			a_{11}(c) &= \frac{1}{N} \ts \Big( (-\gamma_{13} - \gamma_{31})\ts c_1^2 + (\gamma_{12} + \gamma_{21} - \gamma_{13} - \gamma_{31})\ts c_1c_2 + (\gamma_{13} + \gamma_{31} + \gamma_{12}' \\
			&\qquad+ \gamma_{13}' - \gamma_{31}') \ts c_1 + (\gamma_{21}' - \gamma_{31}') \ts c_2 + \gamma_{31}' \Big),
		\end{split}
		\\
		a_{12}(c) &= - \frac{1}{N} \Big( (\gamma_{12} + \gamma_{21}) \ts c_1c_2 + \gamma_{12}'\ts c_1 + \gamma_{21}'\ts c_2  \Big), \\
		\begin{split}
			a_{22}(c) &= \frac{1}{N} \ts \Big( (-\gamma_{23} - \gamma_{32})\ts c_2^2 + (\gamma_{12} + \gamma_{21} - \gamma_{23} - \gamma_{32})\ts c_1c_2 \\
			&\qquad + (\gamma_{23} + \gamma_{32} + \gamma_{21}' + \gamma_{23}' - \gamma_{32}') \ts c_2 + (\gamma_{12}' - \gamma_{32}') \ts c_1 + \gamma_{32}' \Big).
		\end{split}
	\end{align}
\end{subequations}
The remaining entries are given by $a_{13} = a_{11} - a_{12}$, $a_{23} = a_{22} - a_{12}$ and $a_{33} = a_{11} + a_{22} + 2\ts a_{12}$. Using the coefficients appearing in \eqref{eq:red_drift} and \eqref{eq:red_diffusion}, and exploiting
\begin{equation*}
	a_{ij}(x) \approx (\mathcal{L}\psi_k)(x) - b_i(x) x_j - b_j(x) x_i,
\end{equation*}
we can construct the entries of matrix $L_N$. E.g., we obtain $l_{ij} = \gamma_{13} - \gamma_{31}$. We obtain the first columns of $L_{10}$:
\begin{equation*}
	L_{10} =
	\kbordermatrix{
		& 1 & c_1  & c_2    & c_1^2 & c_1c_2 & c_2^2 & \dots \\
		1        & 0 & 0.01 & 0.01   & 0.001 & 0     & 0.001 & \dots  \\
		c_1      & 0 & 0.97 & 0      & 0.321 & 0.009 & 0     & \dots  \\
		c_2      & 0 & 0    & -1.03  & 0     & 0.009 & 0.321 & \dots  \\
		c_1^2    & 0 & -1   & 0      & 1.64  & 0     & 0     & \dots  \\
		c_1c_2   & 0 & -2   & 2      & 0     & -0.36 & 0     & \dots  \\
		c_2^2    & 0 & 0    & 1      & 0     & 0     & -2.36 & \dots  \\
		c_1^3    & 0 & 0    & 0      & -2    & 0     & 0     & \dots  \\
		c_1^2c_2 & 0 & 0    & 0      & -4    & 1     & 0     & \dots  \\
		c_1c_2^2 & 0 & 0    & 0      & 0     & -1    & 4     & \dots  \\
		c_2^3    & 0 & 0    & 0      & 0     & 0     & 2     & \dots
	} \in \mathbb{R}^{10 \times 10}.
\end{equation*}
Note that the indices $ij$ depend on the ordering of the basis elements. Here, it holds that $l_{42} = \gamma_{13} - \gamma_{31}$.

\section{Pseudocode for the discrete-time extended voter model}
\label{apdx:DTEVM}

\begin{algorithm}
	\ForAll{timesteps}{
		\ForAll(randomly){agents}{
			Get number $N$ of adjacent neighbors.\\
			Get number $X_j$ of type $S_j$ in neighborhood for all $d$ types.\\
			Calculate transition probabilities $P = \exp(t_\text{step} \ts G)$ based on neighbors for
			\begin{equation*}
				(G_{ij})_{i, j=1,\dots,d} =
				\begin{cases}
					-\sum_{j=1}^{d} \frac{\gamma_{ij} X_j}{N} + \gamma_{ij}', &\text{if} \quad i=j,\\
					\frac{\gamma_{ij} X_j}{N} + \gamma_{ij}', &\text{else}.
				\end{cases}
			\end{equation*}\\
			Update agent's state according to previously calculated probabilities.
		}
	}
	\caption{Discrete-time extended voter model}
    \label{algo:ABM}
\end{algorithm}

\end{document}